\PassOptionsToPackage{sort&compress}{natbib}
\documentclass[preprint,11pt]{elsarticle}
\usepackage{amsmath,amscd,amsthm,amsxtra,amsfonts,amssymb,algorithm,algorithmicx,algpseudocode,array,arydshln,arydshln,bbm,bm,booktabs,cases,caption,color,enumerate,extarrows,float,graphicx,geometry,hyperref,listings,mathrsfs,manyfoot,multirow,subfig,textcomp,verbatim,xcolor}
\usepackage[utf8]{inputenc}

\allowdisplaybreaks

\usepackage[title]{appendix}%
\numberwithin{equation}{section}

\newtheorem{theorem}{Theorem}
\newtheorem{lemma}[theorem]{Lemma}
\newtheorem*{lemma*}{Lemma}

\theoremstyle{definition}

\newtheorem{remark}{Remark}

\journal{Journal of Computational and Nonlinear Dynamics}

\begin{document}

\begin{frontmatter}

\title{Memory-Dependent FPK Equations for Nonlinear SDOF Oscillators Under Fractional Gaussian Noise Excitation}

\author[1]{Lifang Feng}\ead{FLF.fenglifang@outlook.com}

\author[1,2]{Bin Pei\corref{cor1}}\ead{binpei@nwpu.edu.cn}	

\author[1,2]{Yong Xu}\ead{hsux3@nwpu.edu.cn}

\affiliation[1]{organization={School of Mathematics and Statistics, Northwestern Polytechnical University},
	city={Xi'an},
	postcode={710072},
	country={China}}
\affiliation[2]{organization={MOE Key Laboratory for Complexity Science in Aerospace, Northwestern Polytechnical University},
	city={Xi'an},
	postcode={710072},
	country={China}}

\cortext[cor1]{Corresponding author at: School of Mathematics and Statistics, Northwestern Polytechnical University, Xi'an, 710072, China.}

\begin{abstract}
This paper investigates the transient probabilistic responses of nonlinear single-degree-of-freedom oscillators subjected to external fractional Gaussian noise (FGN) excitation. Owing to the inherent long-range correlations and memory characteristics of FGN, the resulting response process exhibits non-Markovian properties, rendering the classical Fokker-Planck-Kolmogorov (FPK) equation method inapplicable in its direct form. To overcome this critical challenge, a memory-dependent FPK (memFPK) equation is formulated for two-dimensional nonlinear stochastic systems within the fractional Wick-It\^o-Skorohod integral framework. The derived memFPK equation incorporates mixed second-order derivative terms, as well as time-dependent and state-dependent diffusion coefficients, which inherently capture the long-range correlations and memory effects induced by FGN excitation.
For the numerical solution of the memFPK equation, a discretized local mean treatment is developed to estimate the memory-dependent diffusion coefficients involving conditional expectations. The proposed approach integrates local statistical averaging and smoothing techniques to enhance the stability of coefficient estimation. Subsequently, the memFPK equation is numerically solved using a finite difference scheme.
The accuracy and effectiveness of the proposed framework are validated through linear and nonlinear numerical examples. Comparative results demonstrate the excellent agreement with analytical solutions or Monte Carlo simulations in terms of transient joint probability density functions (PDFs), marginal PDFs, low-probability tail regions, and statistical moments. These findings confirm that the proposed memFPK equation method serves as a robust and effective tool for analyzing the transient non-Markovian probabilistic responses of nonlinear SDOF systems under FGN excitation.
\end{abstract}


\end{frontmatter}
\begin{sloppypar}
\section{Introduction}\label{sec1}
Over the past several decades, nonlinear stochastic dynamics has attracted significant research attention due to its widespread applications across engineering and applied sciences. To quantify uncertainty propagation in dynamical systems under random excitations, various analytical and numerical approaches have been developed. Among these, the Fokker-Planck-Kolmogorov (FPK) equation approach \cite{risken1989fokker} stands out as one of the most rigorous methods for characterizing the evolution of probability density functions (PDFs), as it fully captures all statistical properties of system responses. Alternative approaches include: Monte Carlo simulation (MCS) \cite{pradlwarter1997advanced}, the generalized cell mapping method \cite{yue2019Probabilistic,kong2016on,yang2023global}, moment equation techniques \cite{lutes2000direct,di1992stochastic}, statistical linearization methods \cite{spanos2019formulation,kong2021stochastic}, equivalent linearization schemes \cite{proppe2003equivalent,lin2015responses}, stochastic averaging methods \cite{zhu2025stochastic,li2024data}, generalized density evolution equation methods \cite{chen2025efficient,lyu2024refined}, and path integral formulations \cite{psaros2018wiener,mavromatis2024extrapolation,zhu2021stationary,Yai2022efficient,zhao2024an}.

However, many real-world systems, such as those arising in materials, biological processes, and climate systems, exhibit long-range correlations and memory effects. Typical memory-dependent stochastic models incorporate memory characteristics by adopting fractional operators or integro-differential formulations. Relevant studies have systematically examined these systems in terms of theory and numerical implementation \cite{badawi2024existence,badawi2025existence,badawi2025theoretical,badawi2026theoretical}. In engineering stochastic dynamics, fractional-order modeling has also been used to describe memory-dependent material or damping effects. Xu et al. \cite{xu2016method} developed an analytical method for strongly nonlinear stochastic systems with fractional damping, while Liu et al. \cite{liu2018active,liu2020bistability} studied vibration suppression and random fluctuation behavior in airfoil systems with fractional or viscoelastic memory effects.
In the present work, we focus on a different but complementary source of the  long-range correlations  carried directly by the external stochastic excitation.
These excitations can be modeled using fractional Gaussian noise (FGN) \cite{mandelbrot1968fractional}, which is a stationary Gaussian process with  long-range correlation features. FGN has been used in many scientific and engineering fields \cite{vilk2022unravelling,januvsonis2020serotonergic,delgado2007reflected,hristopulos2003permissibility,karacan2009elastic,tey2025distinguishing,kusumayudha2022mathematical,dung2012on,feng2023deep}. Since FGN is non-Markovian, it creates great challenges for analyzing the probabilistic responses of dynamical systems. Although some progress has been made, it only applies to linear systems such as exact solutions for multidimensional linear systems excited by FGN \cite{vyoral2005kolmogorov}, and extensions to nonlinear systems are still limited. For nonlinear systems, Deng and Zhu \cite{lin2015responses} used equivalent linearization to study single-degree-of-freedom (SDOF) nonlinear oscillators under FGN excitation.  Other related progress includes stochastic averaging methods for quasi-Hamiltonian systems under FGN \cite{lu2017stationary,deng2018stochastic}. These methods create approximate stationary PDF by running numerical simulations on averaged stochastic differential equations (SDEs). They focus only on stationary solutions and provide little information about transient behavior. In addition, the classical framework of the FPK equation is based on Markovian assumptions. It cannot be used directly for systems driven by FGN because It\^o calculus does not work in non-Markovian situations.

Recent studies have started investigating the framework of the FPK equation for systems under FGN excitation. Some progress has been reported in \cite{vaskovskii2022analog,choi2021entropy,di2022fokker}, but these works often do not consider drift terms. In order to involve drift terms, Pei et al. \cite{pei2025non} have yielded the memory-dependent probability density evolution equations for nonlinear systems driven by combined FGN and Gaussian white noise (GWN) under the commutativity condition. But this still limits how they can be used in real engineering dynamics. To fix this, we have extended this framework to general conditions without requiring commutativity, thereby establishing a novel memFPK equation that governs the transient PDF of one-dimensional nonlinear systems subjected to FGN and GWN \cite{feng2025memory}. The only minor drawback is that the drift and diffusion terms in the memFPK equation depend on the past conditions and involve conditional expectations. This creates new challenges for numerical computation. Fortunately, the Volterra adjustable decoupling approximation (VADA) method \cite{mamis2019systematic} can be used to estimate these terms. However, our findings show that VADA may not perform well in some nonlinear cases, as it eliminates important memory effects. Meanwhile, it should be specifically noted that the method in \cite{feng2025memory} only applies to the case of one-dimensional nonlinear systems under FGN excitation. It is ineffective for general two-dimensional systems and even SDOF oscillators. 

To address the above challenges, this paper develops a memFPK framework for analyzing the transient probabilistic responses of nonlinear SDOF oscillators under FGN excitation. Since FGN is non-Markovian, the system response does not satisfy the Markovian property required by the Chapman-Kolmogorov equation, so the classical FPK equation cannot be directly used.
Although FGN-driven systems can be transformed into infinite-dimensional Markovian systems using GWN-driven filters \cite{muravlev2011representation}, practical calculations require finite-dimensional truncation \cite{carmona2000approximation}. Low-dimensional truncation often leads to insufficient accuracy, while accurate approximation requires a large number of auxiliary variables and results in high-dimensional augmented systems \cite{bayer2023markovian}. Consequently, the corresponding FPK equation becomes extremely high-dimensional, making it nearly impossible to directly calculate the displacement-velocity joint PDF.
Based on the fractional Wick-It\^o-Skorohod (FWIS) integral theory \cite{ducan2000stochastic,biagini2008stochastic}, this work establishes a new memFPK equation for two-dimensional nonlinear stochastic systems under FGN excitation. This method avoids state-space augmentation, maintains the original system dimension, and allows direct computation of the transient joint PDF in the original displacement-velocity phase space.
	
Notably, while the memFPK equation can realize direct probabilistic evolution without expanding the state space, it still faces great difficulties in numerical calculation.
On this basis, this paper proposes a reliable numerical method to calculate memory-dependent diffusion coefficients, which is the second key contribution of this study. These coefficients include hidden state-related conditional expectations and are easily affected by sparse sampling. Common methods such as VADA often lose key memory information and lead to unstable numerical results in two-dimensional cases. Moreover, if these coefficients do not change smoothly with state variables, regression-based methods will also produce inaccurate extrapolation results.
To solve this problem, we adopt a data-driven discretized local mean treatment, referred to as DLM in this paper. This treatment reconstructs conditional expectations through regional statistical averaging and convolution smoothing. It can steadily estimate memory-dependent diffusion coefficients, and thus supports effective numerical calculation of the memFPK equation for nonlinear SDOF systems.

The proposed memFPK equation distinctly clarifies the distinctions in responses of nonlinear systems excited by GWN and FGN.  For SDOF systems under external GWN excitation, the corresponding FPK equations possess constant diffusion coefficients and exclude mixed second-order partial derivative terms $\partial^2/\partial y_i\partial y_j$. This property is well consistent with the inherent Markovian characteristics of dynamical systems driven by GWN.
On the contrary, when the same systems are subjected to external FGN excitation, the derived memFPK equations show three prominent characteristics. First, mixed second-order derivative terms naturally emerge in the equation formulation. Second, the diffusion coefficients turn to be explicitly time-varying. Third, the diffusion coefficients of nonlinear systems further rely on system state variables.
Such essential differences originate from the long-range correlation induced by FGN. They effectively separate the non-Markovian dynamical behaviors under FGN excitation from the classic Markovian dynamical features governed by GWN. The above theoretical improvements also make the current research fundamentally different from the traditional analytical framework based on conventional Markovian FPK equation methods.

The remainder of the paper is structured as follows. Section 2 outlines the main results. Section 3 explains the proposed DLM treatment in detail, which is used to estimate memory-dependent diffusion coefficients. Section 4 provides several numerical examples to test the accuracy and effectiveness of the method. Finally, Section 5 includes the discussions and conclusions. Appendix introduces the mathematical foundation of the memFPK equation, along with how it is derived using the fractional It\^o's formula and integration by parts.

\section{The memFPK equation for SDOF oscillators excited by FGN}
\subsection{Linear SDOF oscillators}
Consider linear SDOF oscillator excited by FGN, its dynamics is governed as follows
\begin{equation}\label{sdof-linear}
	\ddot{X}(t)+c\dot{X}(t)+kX(t)=\sigma\xi^H(t),
\end{equation}
where $ X,\dot{X},\ddot{X}\in \mathbb{R}$ represent the displacement, velocity, and acceleration, respectively; and a dot over a variable denotes differentiation with respect to time $t$; $c, k$ represent damping and stiffness coefficients, respectively, $\sigma$ is the input force influence; the initial conditions $X(0)=x_0$ and $\dot{X}(0)=\dot{x}_0$ are random initial displacement and velocity variables, respectively; 
$\xi^{H}$ is a unit FGN excitation process with a power-law autocorrelation function, i.e. $$\mathbb{E}[\xi^{H}(t)]=0,\quad \mathbb{E}[\xi^{H}(t)\xi^{H}(t+\tau)]=H(2H-1) \vert\tau\vert ^{2H-2}+2H \vert\tau\vert ^{2H-1}\delta(\tau),$$
in which $\delta(\cdot)$ is Dirac’s delta function, $\mathbb{E}[\cdot]$ is the expectation operator, and $\xi^{H}(t)$ reduces to be GWN when $ H = 1/2 $.

It is pointed out that when $ 0<H<1/2 $, $\xi^{H}(t)$ is not suitable for modeling physical noise since its correlation function is negative. In the present paper, only the case where $ 1/2<H<1 $ is considered. For this range of $H$ (that is, $1/2 < H < 1$), the power spectral density (PSD) of $\xi^{H}(t)$ is 

\begin{align}\label{psd-fgn}
	S_H(\omega)=\frac{H\Gamma(2H)\sin(H\pi)}{\pi} \vert \omega \vert ^{1-2H},
\end{align}
where $ \Gamma(\cdot) $ is the Gamma function, for $ H = 1/2 $, the PSD reduces to a constant power spectrum value equal to $1/(2\pi)$. From Eq.~(\ref{psd-fgn}), we know that $\xi^{H}(t)$ has a power-law PSD at all frequencies $ \omega $ and exhibits long-range correlations and non-Markovian characteristics. 

Set the state vector $\boldsymbol{X}=(X,\dot{X})^{T}=(X, V)^T$, Eq.~(\ref{sdof-linear}) can be expressed in vector form, 
\begin{align*}
	\dot{\boldsymbol X} (t) = A\boldsymbol{X}(t)+ \Sigma \xi^{H}(t),
\end{align*}
where
\begin{align*}
	A=\begin{pmatrix} 0&1 \\ -k&-c\end{pmatrix},
	\Sigma=\begin{pmatrix}  0\\ \sigma \end{pmatrix}.
\end{align*}

For linear SDOF oscillator subjected to FGN excitation, the system responses remain Gaussian. This allows the joint PDF to be fully characterized by its mean vector and covariance matrix; in other words, an analytical expression for the PDF of the solutions to linear SDOF oscillator under FGN excitation is obtainable. 
Denote $\boldsymbol{x}=(x,v)$, thus,
\begin{align}\label{Ex1-exact}
	p(\boldsymbol{x},t)=\frac{1}{2\pi |\Sigma_{\boldsymbol{x}}(t)|^{1/2}}\exp\Bigl\{-\frac{1}{2}(\boldsymbol{x}-\boldsymbol{\mu}_{\boldsymbol{x}}(t))^T\Sigma^{-1}_{\boldsymbol{x}}(t)(\boldsymbol{x}-\boldsymbol{\mu}_{\boldsymbol{x}}(t))\Bigr\},
\end{align}
be the density of the two-dimensional Gaussian distribution $N(\boldsymbol{\mu}_{\boldsymbol{x}}(t),\Sigma_{\boldsymbol{x}}(t))$, with the transient mean vector $\boldsymbol{\mu}_{\boldsymbol{x}}(t)$ and covariance matrix $\Sigma_{\boldsymbol{x}}(t)$ are determined by
\begin{align}\label{Ex1-exact-coff}
	\boldsymbol{\mu}_{\boldsymbol{x}}(t)=&e^{At}\mathbb{E}[\boldsymbol{x}_0]=e^{At}\boldsymbol{\mu}_{0},\cr
	\Sigma_{\boldsymbol{x}}(t)
	=&e^{At}\Sigma_0e^{A^Tt}+H(2H-1)\int_{0}^{t}\int_{0}^{t}e^{A(t-u)}\Sigma\Sigma^{T}e^{A^T(t-v)}\vert u- v\vert ^{2H-2}\mathrm{d}u\mathrm{d}v,
\end{align}
where the initial condition $\boldsymbol{X}(0)=\boldsymbol{x}_0=(x_0,\dot{x}_0)^T$ is a random vector that follows the joint Gaussian distribution $N(\boldsymbol{\mu}_{0},\Sigma_{0})$.

According to \cite{vyoral2005kolmogorov}, the joint response PDF of $X(t)$ and $V(t)$, denoted by $p(x,v,t)$, satisfies the following memFPK equation
\begin{align}\label{memFPK-2d-sdof-linear}
		\frac{\partial }{\partial t}p(x,v,t)=&-\frac{\partial }{\partial x} \big\{ v p(x,v,t)\big\}
		-\frac{\partial }{\partial v}\big\{(-kx-cv)p(x,v,t)\big\}\cr
		&+ \frac{\partial^2 }{\partial v\partial x}\big\{\sigma B_{1}(t)p(x,v,t)\big\}
		+\frac{\partial^2 }{\partial v^2}\big\{\sigma B_{2}(t)p(x,v,t)\big\},
\end{align}
with
\begin{align}\label{memFPK-2d-sdof-linear-b}
	B_{j}(t)=H(2H-1)\sigma\int_0^t[e^{A(t-s)}]_{j2}\vert t-s \vert ^{2H-2} \mathrm{d}s,j=1,2,
\end{align}
where $[e^{A(t-s)}]_{ij}$ denotes its $(i,j)$-th component for $i,j=1,2$.

\subsection{Nonlinear SDOF oscillators}
Consider a nonlinear SDOF oscillator subject to FGN whose motion is governed by the following equation
\begin{equation}\label{SDS}
	\ddot{X}(t)+f(X(t),\dot{X}(t))=\sigma\xi^H(t),
\end{equation}
where $f(X,\dot{X})$ denotes the nonlinear functions of $X,\dot{X}$. Now, we set the state vector $\boldsymbol{X}=(X,\dot{X})^{T}=(X, V)^T$,  $
\boldsymbol{f}(\boldsymbol{X})=(V,-f(X,V))^{T}$, $\Sigma=(0,\sigma)^T$, Eq.~(\ref{SDS}) can be expressed in vector form
\begin{align*}
	\dot{\boldsymbol X} (t) = \boldsymbol{f}(\boldsymbol{X}(t))+ \Sigma  \xi^{H}(t).
\end{align*}

Then, we can obtain the memFPK equation for governing the PDF of Eq.~(\ref{SDS}) as Eq.~(\ref{memFPK-2d-sdof}), for the proofs, readers can see Remark \ref{remark-nonlinear} in Section 3. The joint PDF of $X$ and $V$, satisfies the following memFPK equation
\begin{align}\label{memFPK-2d-sdof}
	\frac{\partial }{\partial t}p(x,v,t)=&-\frac{\partial }{\partial x} \big\{ v p(x,v,t)\big\}-
	\frac{\partial }{\partial v}\big\{-f(x,v)p(x,v,t)\big\}\cr
	&+\frac{\partial^2 }{\partial v\partial x}\big\{\sigma B_{1}(x,v,t)p(x,v,t)\big\}+\frac{\partial^2 }{\partial v^2}\big\{\sigma B_{2}(x,v,t)p(x,v,t)\big\},
\end{align}
with 
\begin{align}\label{memFPK-2d-sdof-b}
	B_{j}(x,v,t)=H(2H-1)\sigma\int_0^t\mathbb{E}[\Psi_{j2}(t;s)\vert X(t)=x,V(t)=v]\vert t-s \vert ^{2H-2} \mathrm{d}s,j=1,2,
\end{align}
where $\Psi(t;s)$ is the state transition matrix satisfies $\dot{\Psi}(t;s)
=\nabla\boldsymbol{f}(\boldsymbol{X}(t))\Psi(t;s),\Psi(s;s)=I$ with the identity matrix, and 
\begin{align*}
	\nabla\boldsymbol{f}=\begin{pmatrix}
		0 & 1 \\
		-\frac{\partial f}{\partial x} & -\frac{\partial f}{\partial v}
	\end{pmatrix}.
\end{align*}

The structure of the FPK equation for SDOF oscillator varies significantly depending on the nature of the stochastic excitation. For instance, consider SDOF oscillator excited by external GWN
\begin{equation}\label{2dsds-gwn}
	\ddot{X}(t)+f(X(t),\dot{X}(t))=\sigma\xi(t),
\end{equation}
where ${\xi}(t)$ is a unit GWN. Then, the PDF of oscillator \eqref{2dsds-gwn}, satisfies the following FPK equation \cite[(6.23)]{sun2006stochastic} 
\begin{align}\label{FPK-2d-sdof-bm1}
	\frac{\partial }{\partial t}p(x,v,t)=-\frac{\partial }{\partial x} \big\{ v p(x,v,t)\big\}-
	\frac{\partial }{\partial v}\big\{-f(x,v)p(x,v,t)\big\}+\frac{\partial^2 }{\partial v^2}\big\{\frac{1}{2}\sigma^2p(x,v,t)\big\}.
\end{align}

\begin{remark}
It can be observed from the FPK equation (\ref{FPK-2d-sdof-bm1}) that for a SDOF oscillator excited by external GWN, its diffusion term only includes the second-order partial derivative with respect to velocity, without any cross-derivative terms. Meanwhile, the corresponding diffusion coefficient keeps constant. In comparison, when the oscillator is subjected to external FGN excitation, the derived memFPK equation (\ref{memFPK-2d-sdof}) possesses three distinct characteristics:
\begin{enumerate}
    \item Nonzero cross-derivative term $\partial^2/(\partial v\,\partial x)$ emerges in the equation;
    \item For linear SDOF oscillators, the memory-related diffusion coefficients are explicitly time-varying via the term $B_j(t)$;
    \item In nonlinear cases, such diffusion coefficients further rely on system state variables, which can be expressed in the form $B_j(x,v,t)$.
\end{enumerate}

Such differences separate the PDF evolution equations of stochastic systems under FGN excitation from those under GWN excitation, which fully reflects the memory effects caused by FGN.
\end{remark}

\subsection{Mathematical support for the memFPK equation}\label{sec3}
Given that linear and nonlinear SDOF oscillators are special cases of following two-dimensional nonlinear system
\begin{align}\label{2-dnonlinear}
		\begin{cases}
			&\dot{Y}_1(t)=f_1(Y_1(t),Y_2(t))+\sigma_{11}\xi_1^{H_1}(t)+\sigma_{12}\xi_2^{H_2}(t),\cr
			&\dot{Y}_2(t)=f_2(Y_1(t),Y_2(t))+\sigma_{21}\xi_1^{H_1}(t)+\sigma_{22}\xi_2^{H_2}(t),
		\end{cases}
\end{align}
where $Y_1,Y_2$ are state variables, $f_1,f_2$ are real-valued nonlinear functions, $\sigma_{11},\sigma_{12},\sigma_{21},\sigma_{22}$ are noise intensities, ${\xi}^{H_1}_1,{\xi}^{H_2}_2$ are two independent unit FGNs. In the following, we will derive the memFPK equation that governs the evolution of the PDF for system (\ref{2-dnonlinear}) under FGN excitation (see Theorem \ref{them-1}). In other words, our results will be applicable to any two-dimensional nonlinear systems with general drift coefficient under external FGN excitation, while the memFPK equation (\ref{memFPK-2d-sdof-linear}) for the linear SDOF oscillator (\ref{sdof-linear}) and the memFPK equation (\ref{memFPK-2d-sdof}) for nonlinear SDOF oscillator (\ref{SDS}) mentioned at the beginning of this section constitute special cases of these findings.

First, we rewrite Eq.~(\ref{2-dnonlinear}) in vector form as
\begin{align}\label{2dsds-1}
	\dot{\boldsymbol Y} (t) = \boldsymbol{\hat f}(\boldsymbol{Y}(t))+ \hat \Sigma \boldsymbol {\xi}^{\boldsymbol H}(t),
\end{align}
where $\boldsymbol{Y}(t)=(Y_1(t),Y_2(t))^T$ is a two-dimensional state variable vector; $\boldsymbol{\hat f}=(f_1,f_2)^{T}$ is a real-valued nonlinear function vector; 
	$$\hat{\Sigma}=
	\begin{pmatrix}
		\sigma_{11} & \sigma_{12}\\
		\sigma_{21} & \sigma_{22}
	\end{pmatrix},$$
	is the noise intensity matrix; $\boldsymbol {\hat \xi}^{\boldsymbol H}=({\xi}^{H_1}_1,{\xi}^{H_2}_2)^{T}$ is the two-dimensional FGN vector including two independent unit FGNs.

Owing to the long-memory property of FGN, the response $\boldsymbol Y(t)$ governed by Eq.~(\ref{2dsds-1}) exhibits non-Markovian behavior. As a result, deriving the joint PDF of $\boldsymbol{Y}(t)=(Y_1(t),Y_2(t))^T$ of system~(\ref{2dsds-1}) using the Chapman-Kolmogorov approach is not feasible and the It\^o calculus framework is inapplicable to FGN. Fortunately, based on the methodologies presented in \cite{pei2025non,feng2025memory}, and by utilizing the fractional It\^o's formula and FWIS integral theory, an approach can be proposed to establish a memFPK equation that governs the transient joint PDF of solution to the two-dimensional nonlinear systems excited by FGN as described in Eq.~(\ref{2dsds-1}). 

We only show the theorem in the main text, with the proofs given in Appendix C. This setup is meant to keep the main text brief and easy to follow, making it simpler for readers to understand the core content quickly.

\begin{theorem}\label{them-1}
	Consider system~(\ref{2dsds-1}) with initial condition $ \boldsymbol{Y}_0=\boldsymbol{y_0} $ which is a random vector with known distribution. Suppose that the function $\boldsymbol{\hat f}$ is differentiable function vector, and $\hat \Sigma$ is a non-zero constant matrix. Then, for the two-dimensional stochastic vector process $\boldsymbol{Y}(t)=(Y_1(t),Y_2(t))^T$, the PDF of $\boldsymbol{Y}(t)$, denoted as $p(y_1,y_2,t)$ or $p(\boldsymbol{y},t)$ with $\boldsymbol{y}=(y_1,y_2)$, satisfies the following memFPK equation
	\begin{align}\label{memFPK-2d-nonlinear-1}
		\frac{\partial }{\partial t}p(\boldsymbol{y},t)=-\sum_{i=1}^{2}\frac{\partial }{\partial y_i}\big\{a_{i}^{{\rm (mem)}}(\boldsymbol{y})p(\boldsymbol{y},t)\big\}+\sum_{i=1}^{2}\sum_{j=1}^{2}\frac{\partial^2 }{\partial y_i\partial y_j}\big\{b_{ij}^{{\rm (mem)}}(\boldsymbol{y},t)p(\boldsymbol{y},t)\big\},
	\end{align}
	where the memory-dependent drift and diffusion coefficients are
    \begin{align}\label{2d-coffs-1}
		a_{i}^{{\rm (mem)}}(\boldsymbol{y})&=f_i(\boldsymbol{y}),\quad i=1,2,\cr
		b_{ij}^{{\rm (mem)}}(\boldsymbol{y},t)&=\sum_{k=1}^{2}\sigma_{ik}B_{jk}(\boldsymbol{y},t),\quad i,j=1,2,
	\end{align}
	with
	\begin{align*}
		B_{jk}(\boldsymbol{y},t)=H_k(2H_k-1)\int_0^t\mathbb{E}[\sigma_{1k}\Psi_{j1}(t;s)+\sigma_{2k}\Psi_{j2}(t;s)\vert \boldsymbol{Y}(t)=\boldsymbol{y}]\vert t-s \vert ^{2H_k-2} \mathrm{d}s,
	\end{align*}
	where $\Psi(t;s)$ is the state transition matrix satisfies $\dot{\Psi}(t;s)
	=\nabla\boldsymbol{f}(\boldsymbol{Y}(t))\Psi(t;s),\Psi(s;s)=I$, and $\nabla\boldsymbol{f}=[\partial f_i/\partial y_j]$ is the Jacobian matrix.
\end{theorem}

\begin{remark}\label{remark-linear-1}
	If we consider the linear SDOF oscillator in Eq.~(\ref{sdof-linear}), Eq.~(\ref{memFPK-2d-sdof-linear}) can be recovered from the general memFPK equation by taking $\boldsymbol{Y}=(X,\dot{X})^{T}=(X,V)^{T}, \boldsymbol{y}=\boldsymbol{x}=(x,v)^{T}$, $\sigma_{11}=\sigma_{12}=\sigma_{21}=0, \sigma_{22}=\sigma, \xi_2^{H_2}=\xi^{H}$. The first FGN $\xi^{H_1}$ does not contribute to the system dynamics because its corresponding noise intensities vanish. In addition,
	\begin{align*}
		\boldsymbol{\hat f} (\boldsymbol{Y})= A \boldsymbol{X} =\begin{pmatrix}
			0 & 1 \\
			-k & -c
		\end{pmatrix}\boldsymbol{X}.
	\end{align*}
	Then, the memory-dependent drift coefficients are given by
\begin{align*}
	a_{1}^{{\rm (mem)}}(\boldsymbol{x})&=v, \, \, a_{2}^{{\rm (mem)}}(\boldsymbol{x})=-kx-cv,
\end{align*}
while the memory-dependent diffusion coefficients reduce to
\begin{align*}
	b_{11}^{{\rm (mem)}}(\boldsymbol{x},t)
	&=\sigma_{11}B_{11}(\boldsymbol{y},t)+\sigma_{12}B_{12}(\boldsymbol{y},t)=0,\cr
	b_{12}^{{\rm (mem)}}(\boldsymbol{x},t)
	&=\sigma_{11}B_{21}(\boldsymbol{y},t)+\sigma_{12}B_{22}(\boldsymbol{y},t)=0,\cr
	b_{21}^{{\rm (mem)}}(\boldsymbol{x},t)
	&=\sigma_{21}B_{11}(\boldsymbol{y},t)+\sigma_{22}B_{12}(\boldsymbol{y},t)=\sigma B_{1}(t),\cr
	b_{22}^{{\rm (mem)}}(\boldsymbol{x},t)
	&=\sigma_{21}B_{21}(\boldsymbol{y},t)+\sigma_{22}B_{22}(\boldsymbol{y},t)=\sigma B_{2}(t).
\end{align*}
These coefficients are fully consistent with those appearing in Eq.~(\ref{memFPK-2d-sdof-linear}). Moreover, in the present linear case, the state-transition matrix $\Psi(t;s)$ reduces to the matrix exponential $e^{A(t-s)}$, which satisfies $\dot{\Psi}(t;s)=A\Psi(t;s)$.
\end{remark}

\begin{remark}\label{remark-nonlinear}
	For the nonlinear SDOF oscillator (\ref{SDS}), Eqs.~(\ref{memFPK-2d-sdof}) and (\ref{memFPK-2d-sdof-b}) can be recovered from the general memFPK equation by taking $\boldsymbol{Y}=\boldsymbol{X}=(X,\dot{X})^{T}=(X,V)^{T}, 
	\boldsymbol{\hat f}(\boldsymbol{Y})=\boldsymbol{f}(\boldsymbol{X})=(V,-f(X,V))^{T}$, $\sigma_{11}=\sigma_{12}=\sigma_{21}=0, \sigma_{22}=\sigma, \xi_2^{H_2}=\xi^{H}$, then, the memory-dependent drift and diffusion coefficients reduce to
	\begin{align*}
			a_{1}^{{\rm (mem)}}(\boldsymbol{x})&=v, \, \, a_{2}^{{\rm (mem)}}(\boldsymbol{x})=-f(x,v),\cr
			b_{11}^{{\rm (mem)}}(\boldsymbol{x},t)&=b_{12}^{{\rm (mem)}}(\boldsymbol{x},t)=0,\, \,
			b_{21}^{{\rm (mem)}}(\boldsymbol{x},t)=\sigma B_{1}(x,v,t),\, \,
			b_{22}^{{\rm (mem)}}(\boldsymbol{x},t)= \sigma B_{2}(x,v,t).
		\end{align*}
	The associated state-transition matrix satisfies $\dot{\Psi}(t;s)
	=\nabla\boldsymbol{f}(\boldsymbol{Y}(t))\Psi(t;s),\Psi(s;s)=I$.
\end{remark}

\section{The estimation of the memory-dependent diffusion coefficients}
Before numerically solving the memFPK equation, it is necessary to first approximate the memory-dependent diffusion coefficients across the relevant state space, as these coefficients are expressed in terms of conditional expectations that are generally non-explicit functions of time and state. Previously, the authors employed the VADA method \cite{mamis2019systematic}, which constructs an approximation through an appropriate stochastic Volterra-Taylor functional series expansion around a specific transient response moment, to address the equivalent memory-dependent drift and diffusion coefficients in one-dimensional nonlinear systems excited by combined FGN and GWN \cite{feng2025memory}. However, in the two-dimensional case involving FGN excitation, which exhibits strong path dependence and nonstationarity, the VADA method becomes inapplicable as it tends to cause a significant loss of historical information.

This study is inspired by the work of Lyu and Chen \cite{lyu2021first, lyu2022unified}, who proposed an approach to estimate the equivalent drift coefficient using data from representative MCS. At each grid point, this estimation is performed via locally weighted scatterplot smoothing (LOWESS). Building on this framework, we propose a numerical estimation treatment called the discretized local mean (DLM), which utilizes state space binning and local statistical averaging, followed by convolution-based smoothing, to construct the equivalent memory-dependent diffusion coefficients.

Specifically, for certain time instant $t_h, h = 1,\cdots, N_t$, some representative MCS datas are obtained. For the $q$-th sample, the responses are denoted as
\begin{align*}
	\begin{cases}
		Y_1(t_h)&=y_{1,q}^{(h)}, \\
		Y_2(t_h)&=y_{2,q}^{(h)}, \\
		D_{t_h}^{H_k} (Y_l(t_n))&=d_{k,l,q}^{(h)},\quad k,l=1,2,
	\end{cases}
\end{align*}for $k=1,...,N$, where $N$ is the number of representative MCS, $d_{k,l,q}^{(h)}$ represents the sample value of $D_{t_h}^{H_k} (Y_l(t_n))$. Here, $t_h$ is omitted for convenience below.

Then, the solution domain $\Omega \in \mathbb{R}^2$ for responses $y_1$ and $y_2$ is uniformly divided into $N_1$ and $N_2$ intervals, respectively. Each resulting rectangular bin $I_{ij}$, where $i=0,1,\cdots,N_1$ and $j=0,1,\cdots,N_2$, centered at $(y_{1,i}, y_{2,j})$, is associated with a subset of the sample trajectories whose states at a given time $t_n$ fall into $I_{ij}$, i.e.,
\begin{align*}
	\mathcal{S}_{ij} = \{ q \vert (y_{1,q}, y_{2,q}) \in I_{ij} \}.
\end{align*}

The bin-wise approximation of the equivalent memory-dependent diffusion coefficients $b_{kl}^{{\rm (ca)}}$ is obtained as the empirical mean over samples in $\mathcal{S}_{ij}$
\begin{align*}
	b_{kl,ij}^{{\rm (ca)}}&=b_{kl}^{{\rm (ca)}}(y_{1,i},y_{2,j},t)=\sigma_{k1}B^{{\rm (ca)}}_{l1}(y_{1,i},y_{2,j},t)
	+\sigma_{k2}B^{{\rm (ca)}}_{l2}(y_{1,i},y_{2,j},t)\cr
	&=\sigma_{k1}\mathbb{E}[D_{t}^{H_1}(Y_l(t))\vert Y_1(t)=y_{1,i},Y_2(t)=y_{2,j}]
	+\sigma_{k2}\mathbb{E}[D_{t}^{H_2}(Y_l(t))\vert Y_1(t)=y_{1,i},Y_2(t)=y_{2,j}]\cr
	&= \frac{\sigma_{k1}}{\vert\mathcal{S}_{ij}\vert} \sum_{q \in \mathcal{S}_{ij}} d_{1,l,q}
	+\frac{\sigma_{k2}}{\vert\mathcal{S}_{ij}\vert} \sum_{q \in \mathcal{S}_{ij}} d_{2,l,q},
\end{align*}provided that $\vert\mathcal{S}_{ij}\vert$ is non-zero, and $b_{kl,ij}^{{\rm (ca)}}=\frac{\sigma_{k1}}{N} \sum_{k=1}^{N} d_{1,l,q}+\frac{\sigma_{k2}}{N} \sum_{k=1}^{N} d_{2,l,q}$ when $\vert\mathcal{S}_{ij}\vert = 0$.

To ensure spatial continuity and decrease numerical artifacts caused by local sampling, a two-dimensional smoothing operation is applied to the memory-dependent diffusion coefficients $b_{kl,ij}^{{\rm (ca)}}$. Specifically, a discrete uniform convolution kernel is employed, which is defined as
\begin{align*}
	K(m,n) = \frac{1}{(2r+1)^2}, \quad \text{for } m,n \in {-r, \cdots, r}.
\end{align*}

The smoothed equivalent memory-dependent diffusion coefficients are then given by
\begin{align*}
	\tilde{b}_{kl,ij}^{{\rm (ca)}} = \sum_{m=-r}^{r} \sum_{n=-r}^{r} K(m,n) \cdot b_{kl,i+m,j+n}^{{\rm (ca)}}.
\end{align*}
This smoothing step ensures that $\tilde{b}_{kl,ij}^{{\rm (ca)}}$ vary smoothly across the state space, thereby facilitating interpolation and numerical differentiation in later computational stages.

\begin{remark}
	Here are several points to note at the end of this section:
	\begin{itemize}
		\item To improve computational efficiency, the DLM treatment can be run on a roughly discretized grid in the state space. After that, it can be interpolated onto finer grids using common interpolation methods such as linear or cubic spline interpolation. This approach allows for accurate and efficient reconstruction of the equivalent memory-dependent diffusion coefficients across the entire domain, while effectively reducing excessive computational costs.
		\item The proposed DLM treatment has computational efficiency and clear structural interpretability. It avoids the need for pointwise regression by using bin-wise averaging combined with uniform kernel smoothing. This makes it possible to directly use local statistics obtained from a limited set of deterministic trajectories in a fully non-parametric and model-free way.
		\item However, the proposed treatment currently has certain limitations in its numerical implementation. Specifically, near the domain boundaries, some bins may contain very few or even no samples, which can lead to less accurate estimates in these areas. As a result, the accuracy of the results may be affected in the tails of the response PDF. Possible future improvements could include adding adaptive strategies or deep learning algorithms to enhance the reliability and consistency of estimates when samples are sparse.
	\end{itemize}
\end{remark}

\section{Numerical results}
In the following section, the effectiveness and accuracy of the proposed memFPK equation will be verified through some numerical examples.

\subsection{Example 1}
Consider a linear SDOF oscillator (\ref{sdof-linear}) excited by FGN. The system parameters are chosen as $k=1, c=0.4, \boldsymbol{\mu}_{0}=(-1,-1)^T,\Sigma_{0}={\rm diag}(0.15,0.15)$. The stochastic excitation is modeled by FGN with a Hurst index $H=0.8$, reflecting long-range temporal dependence. Applying \cite[Corollary 2.4]{bernstein1993some}, the explicit expression of $e^{At}$ can be written as
\begin{align}\label{eAt-ana}
	e^{At}=e^{-\zeta\omega_{n}t}
	\begin{pmatrix}
		{\cos(\omega_{d}t)+\frac{\zeta}{\sqrt{1-\zeta^{2}}}\sin(\omega_{d}t)} & {\frac{1}{\omega_{d}}\sin(\omega_{d}t)} \\
		{\frac{-\omega_{d}}{1-\zeta^{2}}\sin(\omega_{d}t)} & {\cos(\omega_{d}t)-\frac{\zeta}{\sqrt{1-\zeta^{2}}}\sin(\omega_{d}t)}
	\end{pmatrix},
\end{align}
where $w_n=\sqrt{k}, \zeta=c/(2\sqrt{k}), w_d=w_n\sqrt{1-\zeta^2}$. Then, substituting Eq.~(\ref{eAt-ana}) into Eqs.~(\ref{Ex1-exact})-(\ref{memFPK-2d-sdof-linear-b}) yields the analytical joint PDF solution via mean vector and covariance matrix and the corresponding memFPK equation, respectively.

\begin{figure}[t!]
	\centering
	\subfloat[memFPK $t=$ 1.0 s]{
		\includegraphics[scale=0.35]{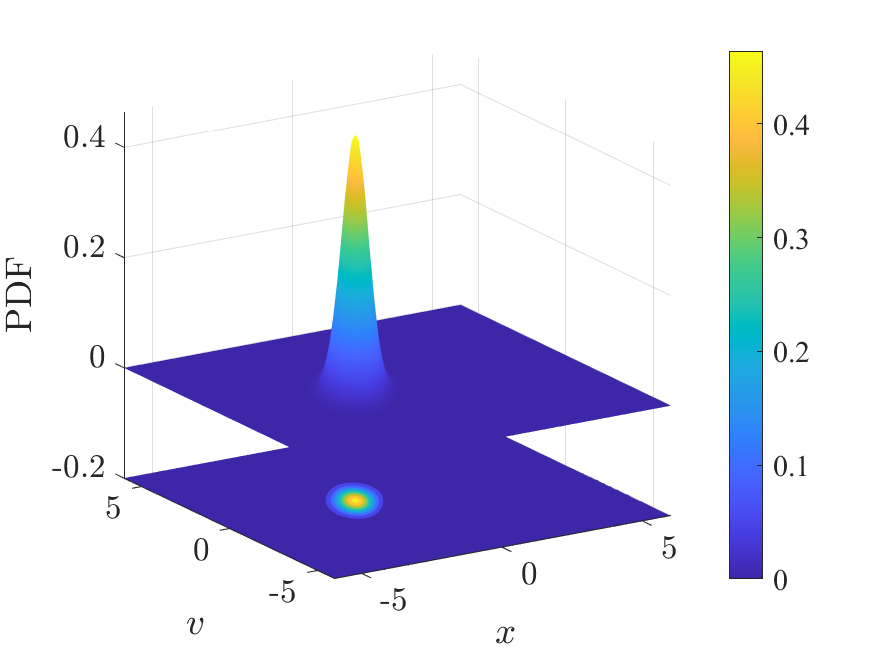}}
	\subfloat[memFPK $t=$ 5.0 s]{
		\includegraphics[scale=0.35]{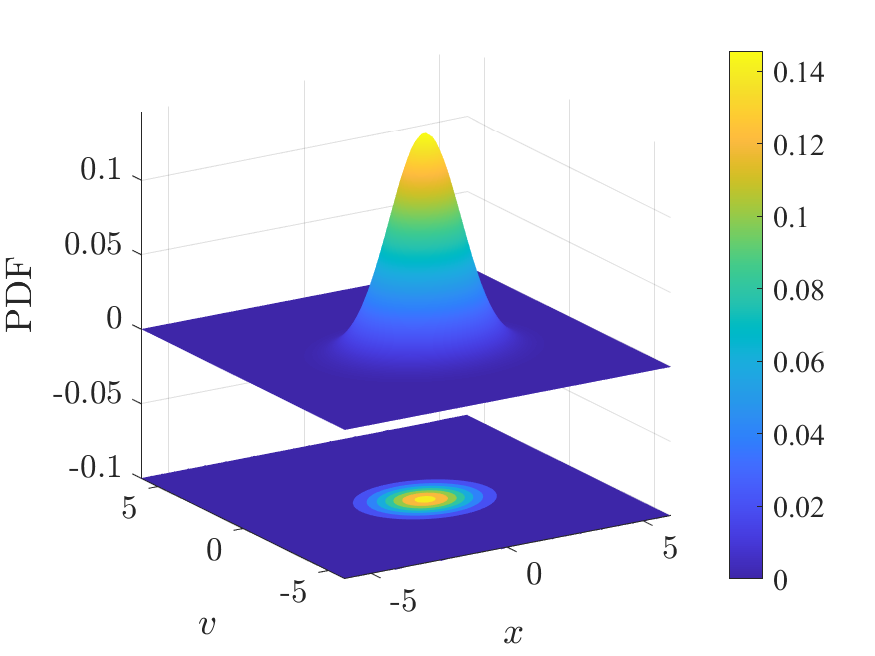}}
	\subfloat[memFPK $t=$ 20.0 s]{
		\includegraphics[scale=0.35]{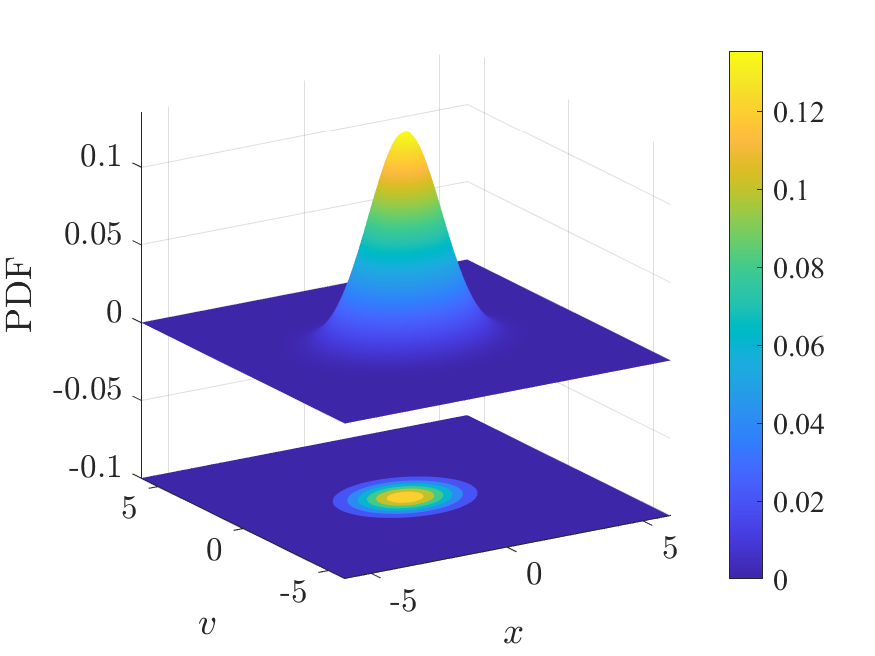}}
	\\
	\subfloat[Analytical $t=$ 1.0 s]{
		\includegraphics[scale=0.35]{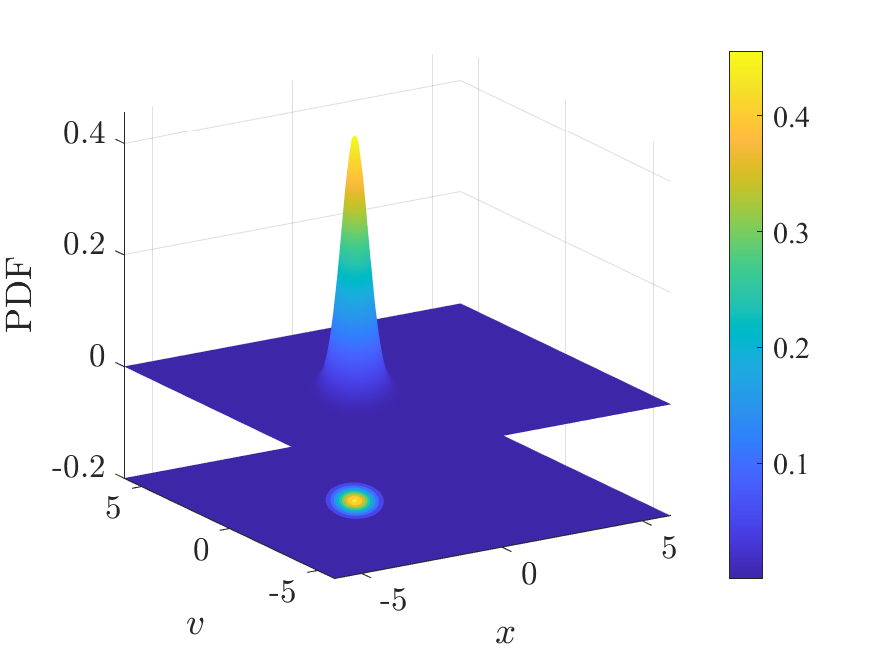}}
	\subfloat[Analytical $t=$ 5.0 s]{
		\includegraphics[scale=0.35]{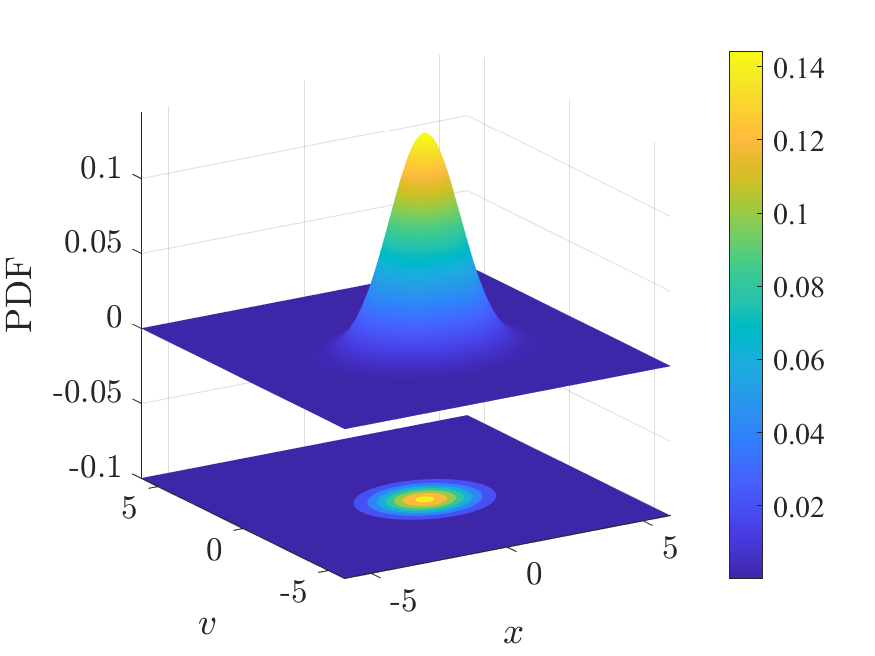}}
	\subfloat[Analytical $t=$ 20.0 s]{
		\includegraphics[scale=0.35]{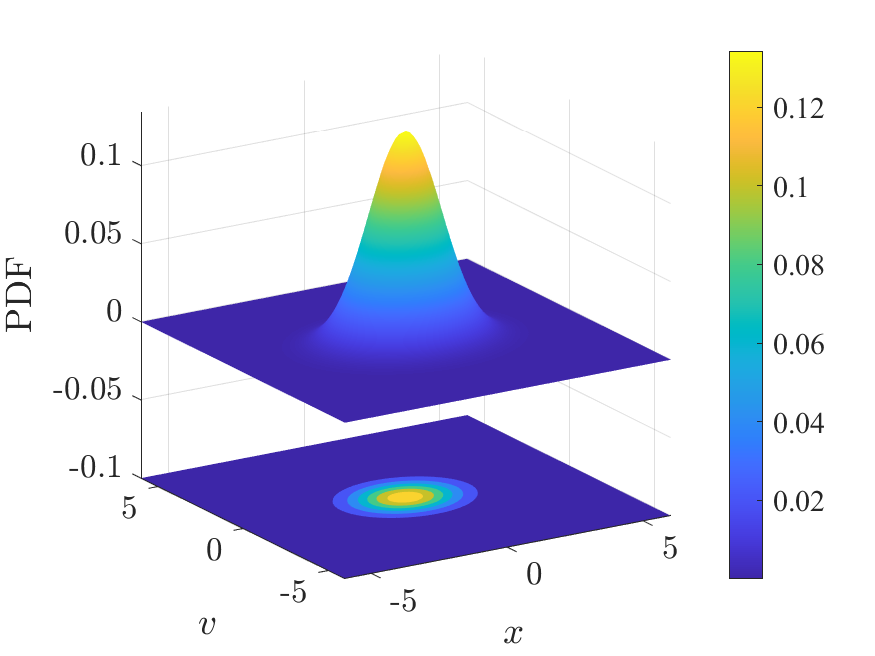}}
	\caption{Joint PDFs of Example 1 obtained from the memFPK equation and the analytical solution obtained from the mean vector and covariance matrix: (a)-(c) memFPK equation solutions; (d)-(f) analytical solutions.}
	\label{fig1-1}
\end{figure}

The covariance matrix in the form of a double integral Eq.~(\ref{Ex1-exact-coff}) is calculated numerically using Wolfram Mathematica 11.3, and this solves the analytical transient joint PDF (\ref{Ex1-exact}). Figs.~\ref{fig1-1}(d)-(f) present the three-dimensional surface plots together with the corresponding contour plots of the analytical joint PDF at representative time instants $t=$ 1.0, 5.0, and 20.0 s.

To assess the accuracy of the memFPK equation, the FD method is employed to obtain its numerical solution. The computational domain is set as $x \in [-6, 6]$ and $v \in [-6, 6]$, with spatial step sizes $\Delta x = 0.15$ and $\Delta v = 0.15$, and time step size $\Delta t = 0.001$ s. Figs.~\ref{fig1-1} (a)-(c) present the numerical joint PDF at the same moment, plotted as both surface and contour plots. It can be observed from Fig.~\ref{fig1-1} that the memFPK equation solutions exhibit excellent agreement with the analytical solutions in terms of peak amplitudes, geometric structures, and spatial decay behaviors.

\begin{figure}[t!]
	\centering
	\subfloat[$t=$ 1.0 s]{
		\includegraphics[scale=0.35]{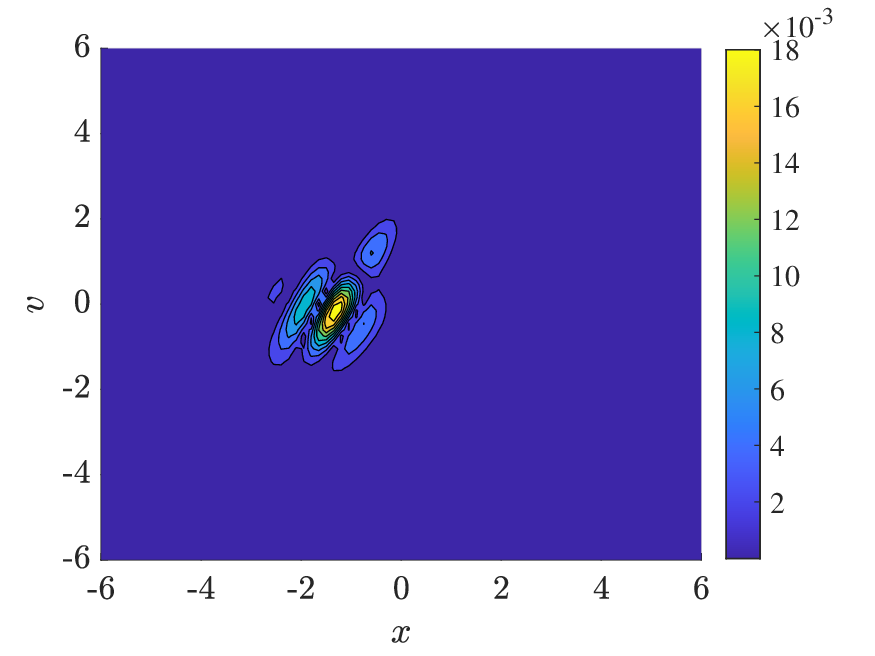}}
	\subfloat[$t=$ 5.0 s]{
		\includegraphics[scale=0.35]{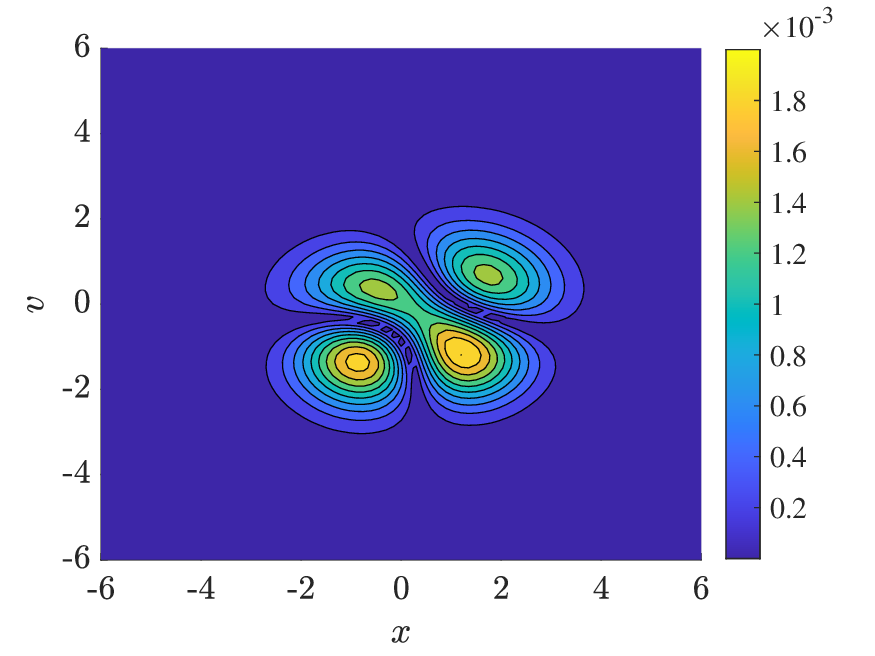}}
	\subfloat[$t=$ 20.0 s]{
		\includegraphics[scale=0.35]{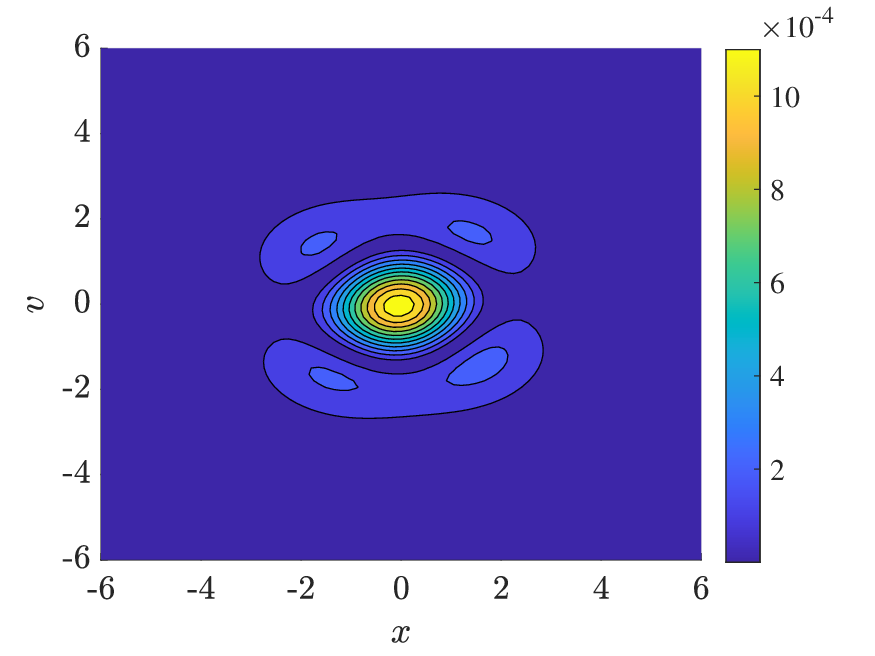}}
	\caption{Absolute error contours of the joint PDF in Example 1 between the memFPK equation solution and the analytical solution obtained from the mean vector and covariance matrix.}
	\label{fig1-2}
\end{figure}

To further quantify the accuracy, Fig.~\ref{fig1-2} shows the absolute error between the joint PDF derived from the memFPK equation and analytical solution via mean vector and covariance matrix. The largest error values mainly appear in the central region where the PDF values attain their maximum magnitudes. However, the overall error remains within the range of $10^{-3}$ to $10^{-4}$,  which indicates that the numerical accuracy is satisfactory throughout the computational domain.

\begin{figure}[t!]
	\centering
	\subfloat[$X(t)$ linear scale]{
		\includegraphics[scale=0.5]{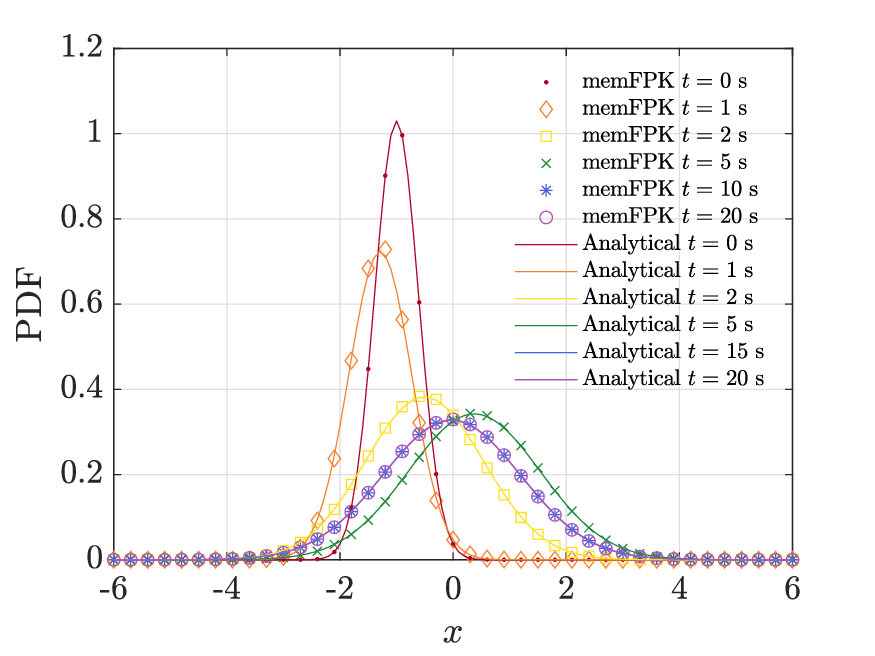}}
	\subfloat[$X(t)$ logarithmic scale]{
		\includegraphics[scale=0.5]{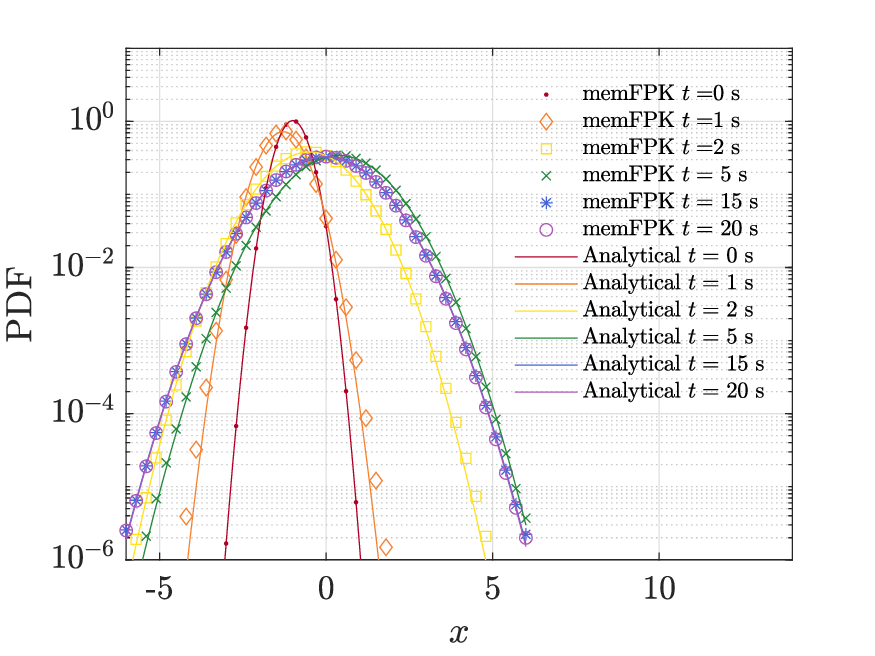}}
	\\
	\subfloat[$V(t)$ linear scale]{
		\includegraphics[scale=0.5]{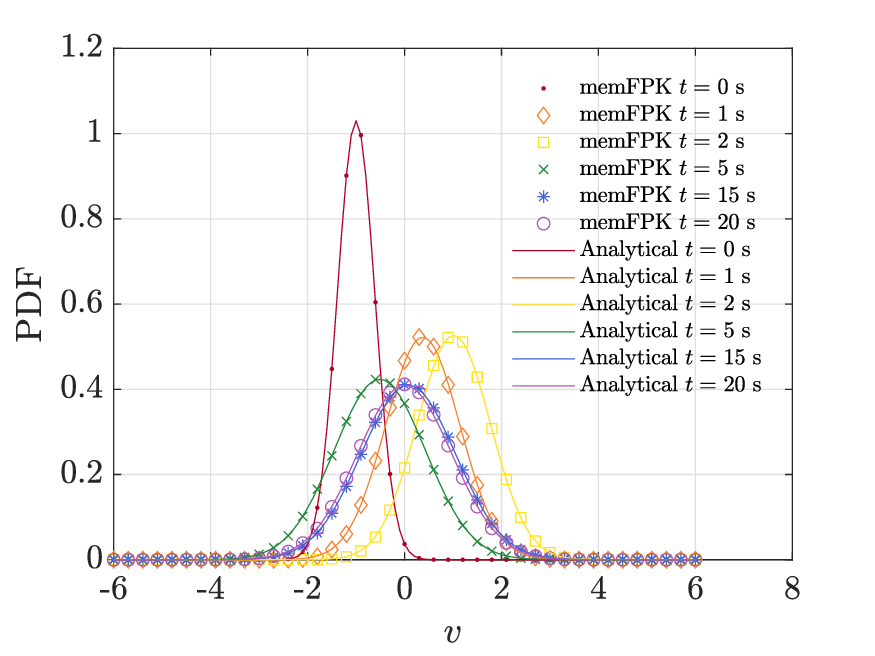}}
	\subfloat[$V(t)$ logarithmic scale]{
		\includegraphics[scale=0.5]{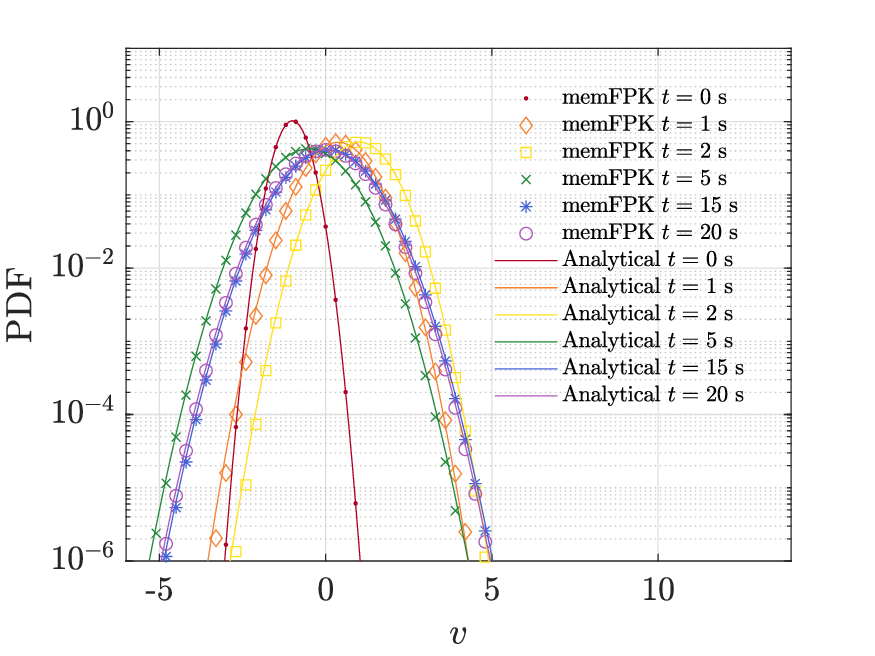}}
	\caption{Transient marginal PDFs of $X(t)$ and $V(t)$ in Example 1 obtained from the memFPK equation solution and analytical solution obtained from the mean vector and covariance matrix.}
	\label{fig1-3}
\end{figure}

Fig.~\ref{fig1-3} presents the transient marginal PDFs of $X(t)$ and $V(t)$. These PDFs are obtained from solving the memFPK equation and compared with the analytical solutions obtained from the mean vector and covariance matrix. Both linear and logarithmic scale representations are employed to examine the agreement in both the high-probability central region and the low-probability tail region. In particular, Figs.~\ref{fig1-3}(b) and (d) demonstrate excellent agreement even in the low-probability tail regions. Since the tail behavior is particularly sensitive to numerical diffusion and approximation errors, the observed agreement further confirms the the accuracy and reliability of the numerical solution to the memFPK equation.

\subsection{Example 2}
Consider a response process $X$ that is governed by the following nonlinear SDOF Duffing oscillator excited by FGN		
\begin{equation}\label{Ex2}
	\ddot{X}+\eta\dot{X}+\alpha X+\beta X^3=\sigma\xi^H(t),
\end{equation}
where $\eta>0$ represents the damping coefficient; $\alpha, \beta$ represent the linear and nonlinear stiffness coefficients, respectively; $\sigma\xi^H(t) $ is a zero-mean FGN term, which models the excitation induced by random perturbations of intensity $\sigma$. The initial condition $\boldsymbol{X}(0)=\boldsymbol{x}_0=(x_0,\dot{x}_0)^T$ is random vector that follows the joint Gaussian distribution $N(\boldsymbol{\mu}_{0},\Sigma_{0})$.

Set the state vector $\boldsymbol{X}=(X,\dot{X})^{T}=(X, V)^T$, Eq.~(\ref{Ex2}) can be expressed in vector form, where
\begin{align*}
	\boldsymbol{f}(\boldsymbol{X})=\begin{pmatrix}V \\ -\eta V-\alpha X-\beta X^3\end{pmatrix},
	\Sigma=\begin{pmatrix}  0\\ \sigma \end{pmatrix}.
\end{align*}

Then, the response joint PDF of oscillator~(\ref{Ex2}) satisfies the memFPK equation \eqref{memFPK-2d-sdof} with the memory-dependent coefficients (\ref{memFPK-2d-sdof-b}), where
$$\nabla\boldsymbol{f}(\boldsymbol{X})=\begin{pmatrix}
	0 & 1 \\
	-\alpha-3\beta X^2 & -\eta
\end{pmatrix}.$$ 
The parameters settings of this oscillator are $\eta = 1, \alpha = -1, \beta = 1, \sigma^2 = 0.36, \boldsymbol{\mu}_{0}=(0,0)^T,\Sigma_{0}={\rm diag}(0.05,0.05), H=0.65$. 

Herein, based on 2000 representative MCS samples, the DLM treatment is employed to construct the equivalent memory-dependent diffusion coefficients at each time step. Subsequently, this coefficients are substituted into the memFPK equation, which is solved numerically using the FD method to obtain the transient joint PDF of $X(t)$ and $V(t)$. Each sample for MCS is simulated by the second order stochastic Runge-Kutta algorithm algorithm \cite{hong2021optimal} for FGN excitation. 

The computational domain is set as $x\in[-2.5, 2.5]$ and $v\in[-2.5, 2.5]$, with spatial step sizes of $\Delta x=\Delta v = 0.083$ and a time step of $\Delta t =$ 0.001 s.

\begin{figure}[t!]
	\centering
	\subfloat[memFPK $t=$ 2.5 s]{
		\includegraphics[scale=0.35]{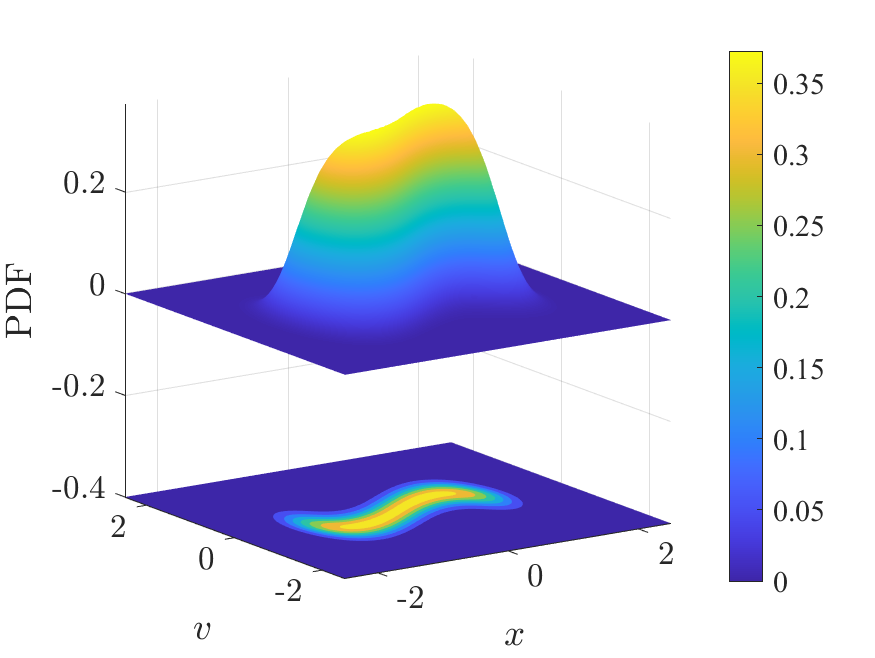}}
	\subfloat[memFPK $t=$ 3.0 s]{
		\includegraphics[scale=0.35]{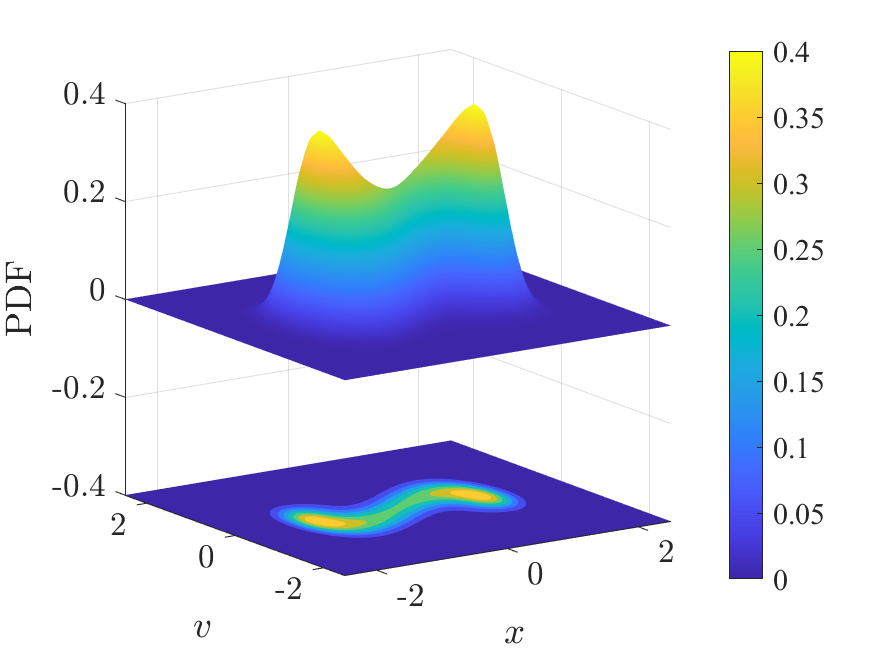}}
	\subfloat[memFPK $t=$ 8.0 s]{
		\includegraphics[scale=0.35]{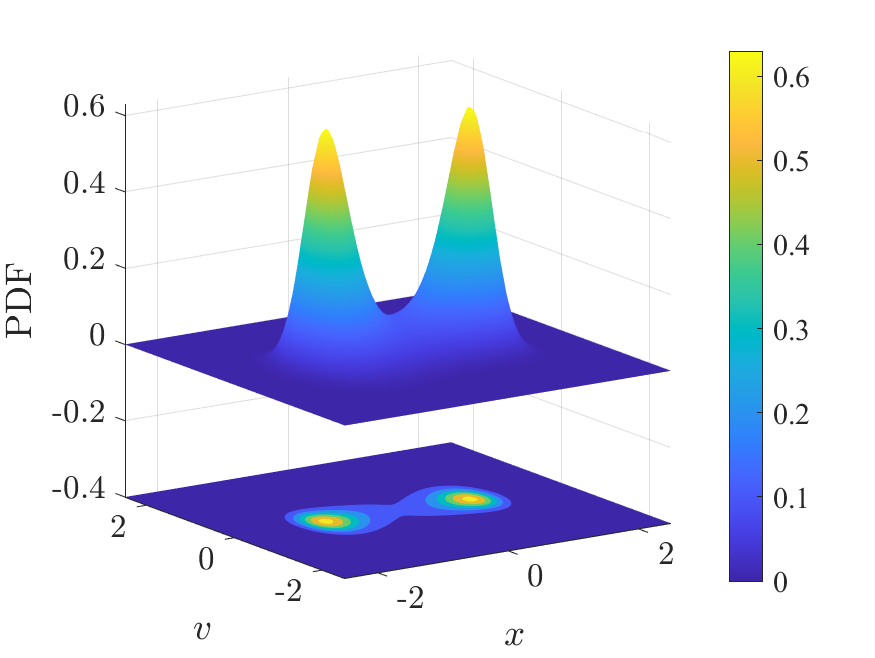}}
	\\
	\subfloat[MCS $t=$ 2.5 s]{
		\includegraphics[scale=0.35]{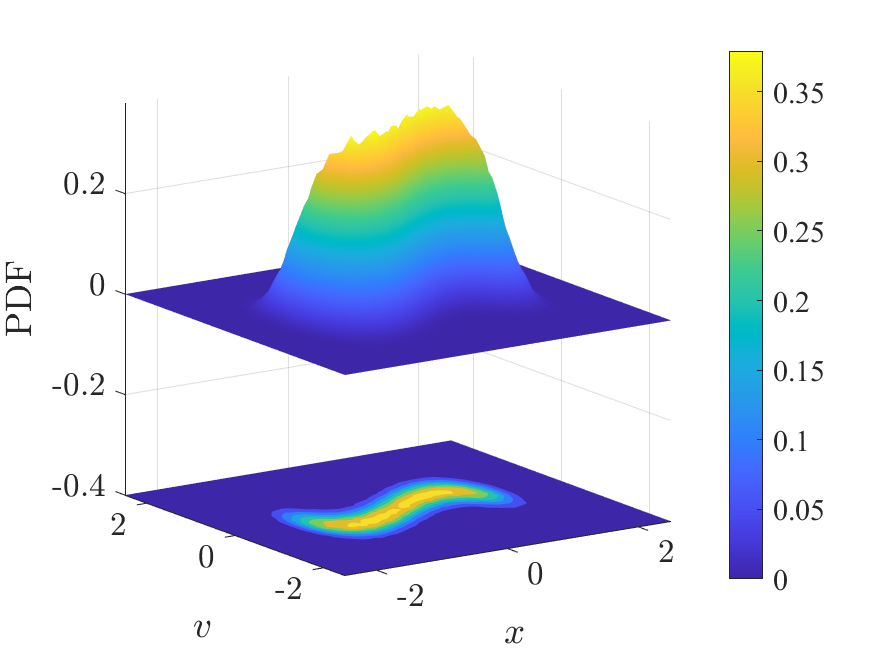}}
	\subfloat[MCS $t=$ 3.0 s]{
		\includegraphics[scale=0.35]{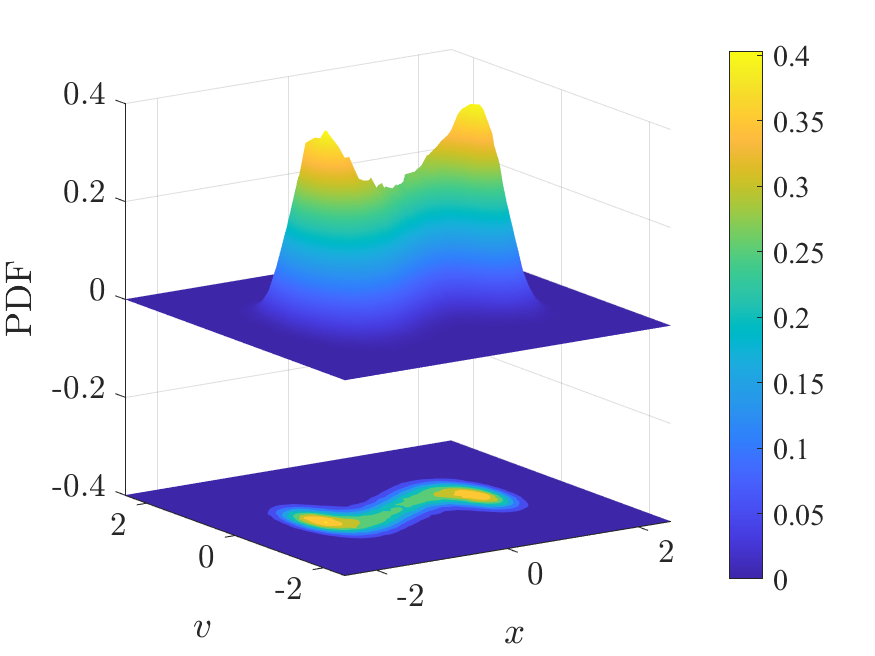}}
	\subfloat[MCS $t=$ 8.0 s]{
		\includegraphics[scale=0.35]{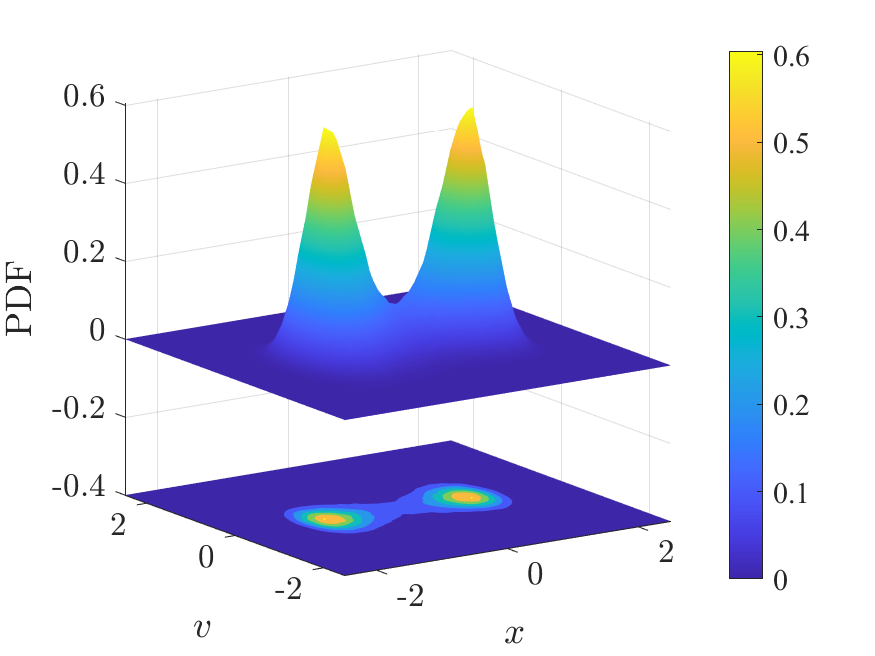}}
	\caption{Joint PDF of Example 2 obtained from memFPK equation and MCS solutions ($ 6\times10^6 $ samples): (a)-(c) memFPK equation solutions vs (d)-(f) MCS solutions.}
	\label{fig2-1}
\end{figure}

\begin{figure}[t!]
	\centering
	\subfloat[$X(t)$ linear scale]{
		\includegraphics[scale=0.5]{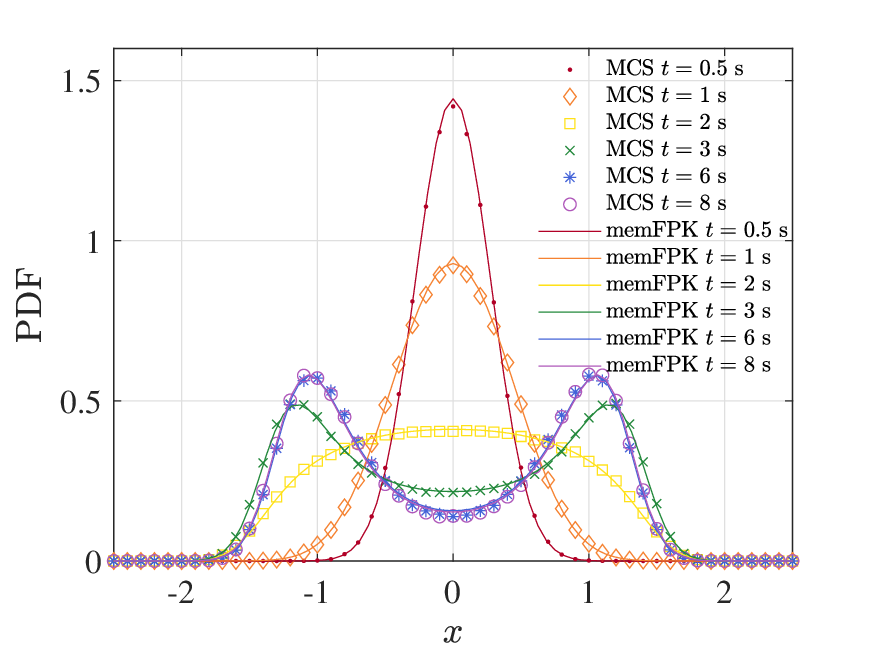}}
	\subfloat[$X(t)$ logarithmic scale]{
		\includegraphics[scale=0.5]{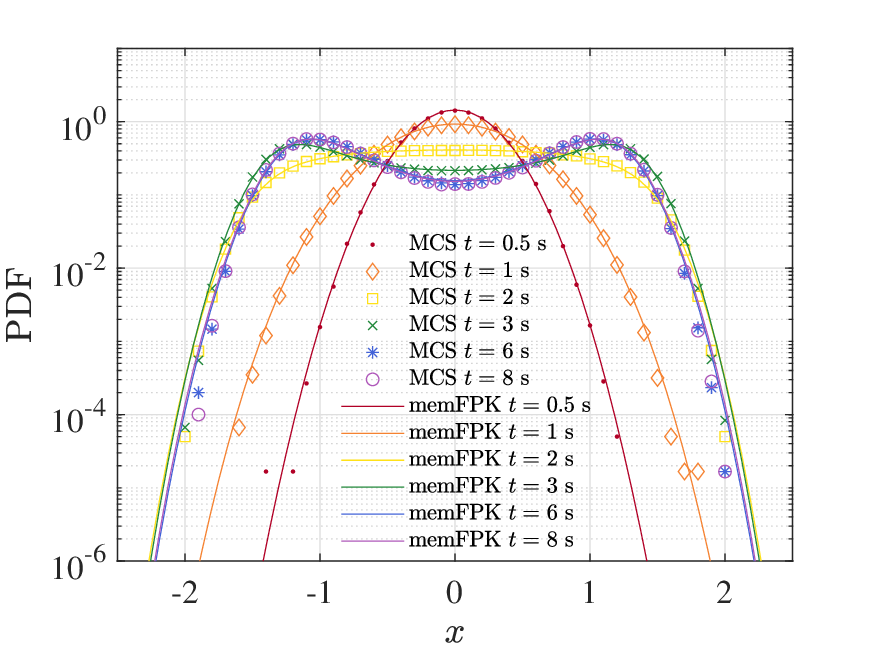}}
	\\
	\subfloat[$V(t)$ linear scale]{
		\includegraphics[scale=0.5]{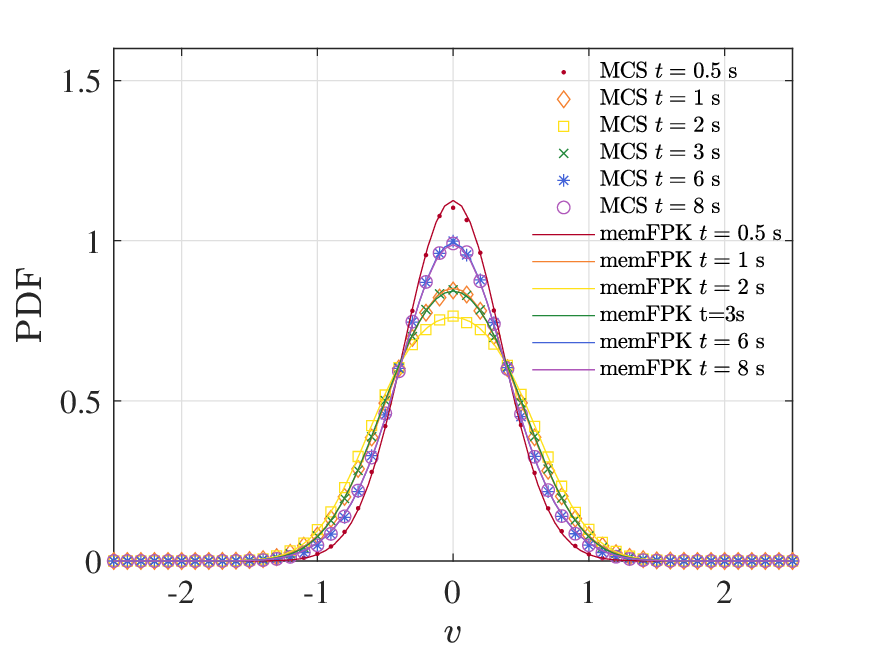}}
	\subfloat[$V(t)$ logarithmic scale]{
		\includegraphics[scale=0.5]{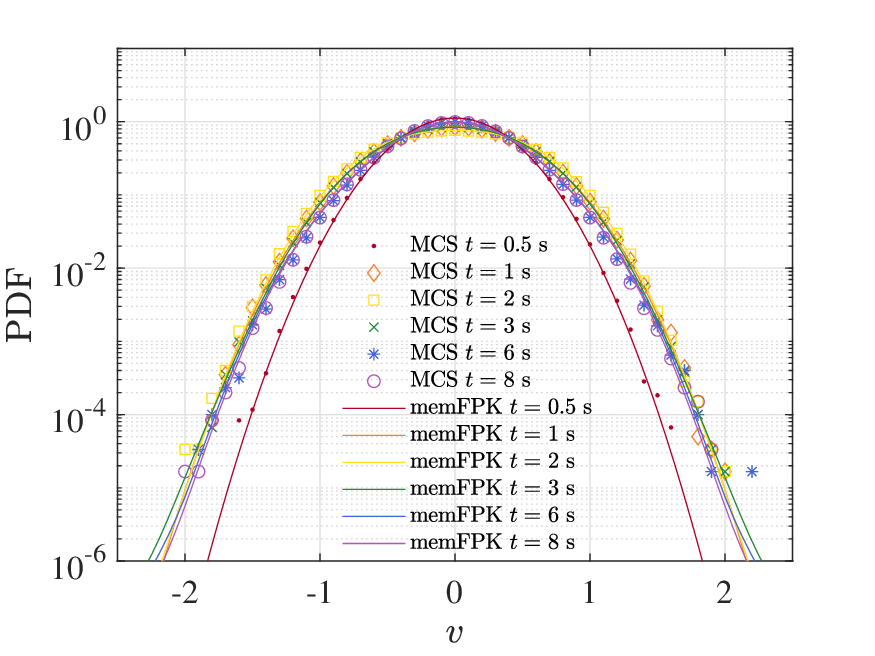}}
	\caption{Transient marginal PDFs of $X(t)$ and $V(t)$ of Example 2 obtained from memFPK and MCS solutions ($ 6\times10^6 $ samples).}
	\label{fig2-3}
\end{figure}

\begin{figure}[t!]
	\centering
	\subfloat[$X(t)$]{
		\includegraphics[scale=0.45]{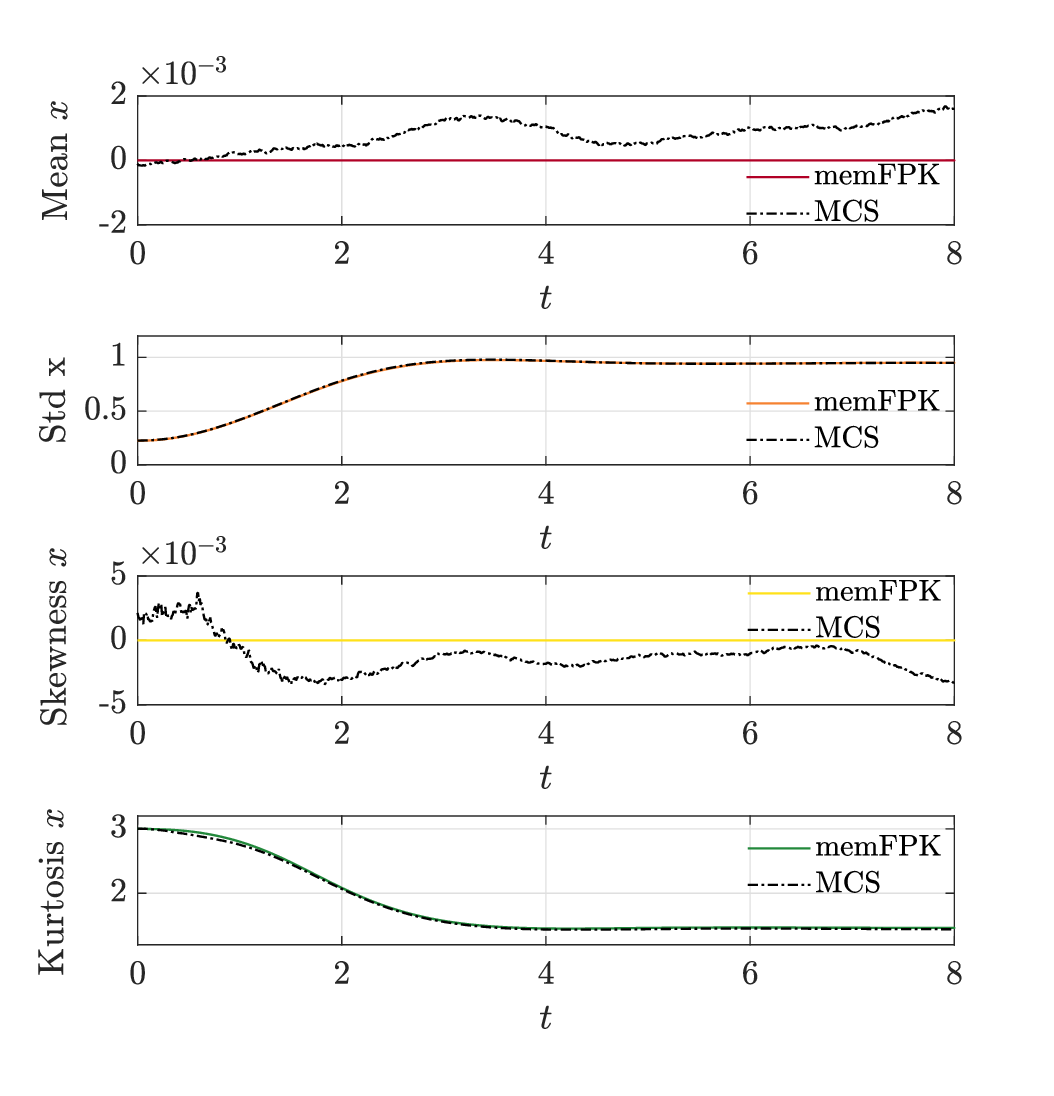}}
	\subfloat[$V(t)$]{
		\includegraphics[scale=0.45]{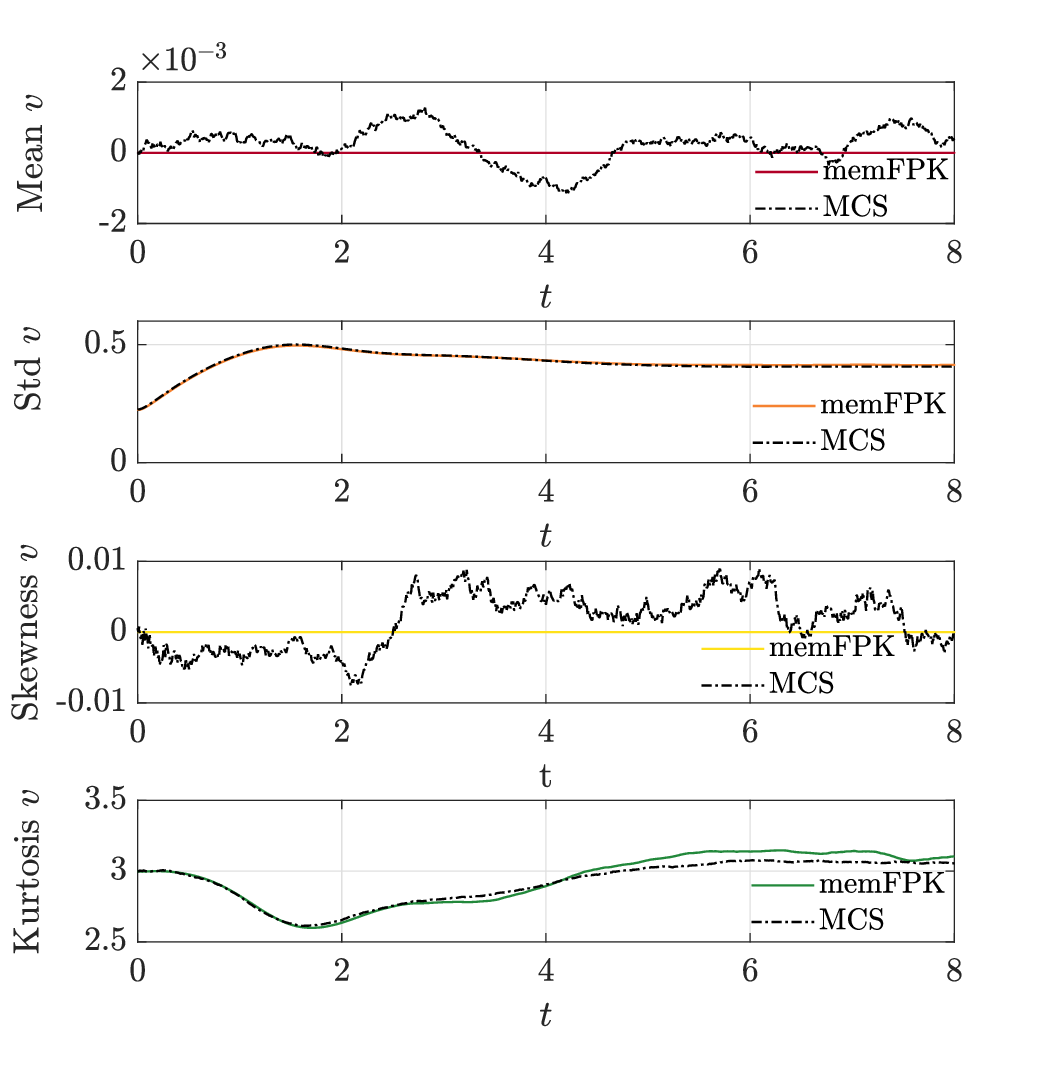}}
	\caption{Temporal evolution of the statistical moments of $X(t)$ and $V(t)$ of Example 2 obtained from memFPK equation and MCS solutions ($ 6\times10^6 $ samples).}
	\label{fig2-5}
\end{figure}

For nonlinear stochastic dynamical system excited by FGN, the exact steady-state or transient solutions are generally unavailable. Therefore, MCS results are used as reference solution in the following examples. The comparisons of the three-dimensional surface and corresponding contour plots of the joint PDF obtained by the memFPK equation method and the MCS method (using \(6\times10^{6}\) samples) at representative time instants $t=$ 2.5, 3.0, and 8.0 s are shown in Fig.~\ref{fig2-1}.  As shown in Fig.~\ref{fig2-3}, the comparisons of the transient marginal PDF of $X(t)$ and $V(t)$ in linear and logarithmic coordinates at typical time points. Fig.~\ref{fig2-5} further presents a comparative analysis of the temporal evolution of the first four statistical moments obtained from both approaches.

The results in Figs.~\ref{fig2-1}-\ref{fig2-5} indicate that the memFPK equation method agrees well with the MCS method in terms of the transient joint PDFs, marginal PDFs, as well as the statistical moments. As shown in Figs.~\ref{fig2-1} and \ref{fig2-3}, the response probability gradually separates into two regions along the displacement direction, and the marginal PDF of $X(t)$ develops a clear bimodal structure during the transient evolution. This behavior reflects the bistable potential structure of the Duffing oscillator under the selected parameter setting, for which probability mass gradually spreads toward the two potential wells under FGN excitation. In contrast, the marginal PDF of $V(t)$ remains essentially unimodal, indicating that the velocity response is still mainly concentrated around zero, even though the displacement response becomes distributed around the two wells. Furthermore, Figs.~\ref{fig2-3}(b) and (d) show that the proposed memFPK equation method maintains excellent agreement with the MCS results even in the low-probability tail regions where the PDF values decrease to the order of $10^{-6}$. By comparison, the MCS results based on $6\times10^{6}$ samples can resolve the tail region only to approximately the order of $10^{-4}$-$10^{-5}$, whereas the proposed method provides more accurate tail descriptions using only 2000 representative samples.

The evolution of the first four statistical moments in Fig.~\ref{fig2-5} provides further information on the transient probabilistic response of the Duffing oscillator. The mean values and skewness of both $X(t)$ and $V(t)$ remain close to zero, indicating the symmetric response distributions. The standard deviation of $X(t)$ increases and then approaches a bounded level as the probability mass spreads toward the two wells. Meanwhile, the kurtosis of $X(t)$ decreases markedly below the Gaussian value of 3, which is consistent with the formation of a bimodal non-Gaussian displacement distribution. These results further demonstrate that the proposed memFPK equation method can accurately capture the transient non-Gaussian features and bimodal probabilistic behavior of the Duffing oscillator.


\subsection{Example 3}
Consider a response process $X$ that is governed by the following nonlinear SDOF Van der Pol oscillator excited by FGN	
\begin{equation}\label{Ex3}
	\ddot{X}+\eta(-1+X^2+\dot{X}^2)\dot{X}+X=\sigma\xi^H(t),
\end{equation}
where $\eta$ represents the control parameters of nonlinearity and damping intensity; $\sigma\xi^H(t)$ represents FGN with the noise intensity $\sigma$. The initial condition $\boldsymbol{X}(0)=\boldsymbol{x}_0=(x_0,\dot{x}_0)^T$ is random vector that follows the joint Gaussian distribution $N(\boldsymbol{\mu}_{0},\Sigma_{0})$.

Set the state vector $\boldsymbol{X}=(X,\dot{X})^{T}=(X, V)^T$, Eq.~(\ref{Ex3}) can be expressed in vector form, where
\begin{align*}
	\boldsymbol{f}(\boldsymbol{X})=\begin{pmatrix}V \\ -\eta(-1+X^2+V^2)V-X\end{pmatrix},
	\Sigma=\begin{pmatrix}  0\\ \sigma \end{pmatrix}.
\end{align*}

Then, the response joint PDF of oscillator~(\ref{Ex3}) satisfies the memFPK equation \eqref{memFPK-2d-sdof} with the memory-dependent coefficients (\ref{memFPK-2d-sdof-b}), where
$$\nabla\boldsymbol{f}(\boldsymbol{X})=\begin{pmatrix}
	0 & 1 \\
	-2\eta XV-1 & -\eta(-1+X^2+3V^2)
\end{pmatrix}.$$ 
The parameter settings of this example are $\eta = 2, \sigma^2 = 1, \boldsymbol{\mu}_{0}=(0,0)^T,\Sigma_{0}={\rm diag}(0.05,0.05), H=0.6$. Similar to Example 2, the DLM treatment is applied to estimate the memory-dependent diffusion coefficients using 2000 representative MCS samples. The memFPK equation with the estimated diffusion coefficients substituted in is solved by FD method. The computational domains are set as $x \in [-2.5, 2.5]$ and $v \in [-2.5, 2.5]$, with spatial step sizes of $\Delta x = 0.0625$ and $\Delta v = 0.0625$, respectively. The time step size $\Delta t = $ 0.001 s is used in FD for time discretization.

\begin{figure}[t!]
	\centering
	\subfloat[memFPK $t=1.0$ s]{
		\includegraphics[scale=0.35]{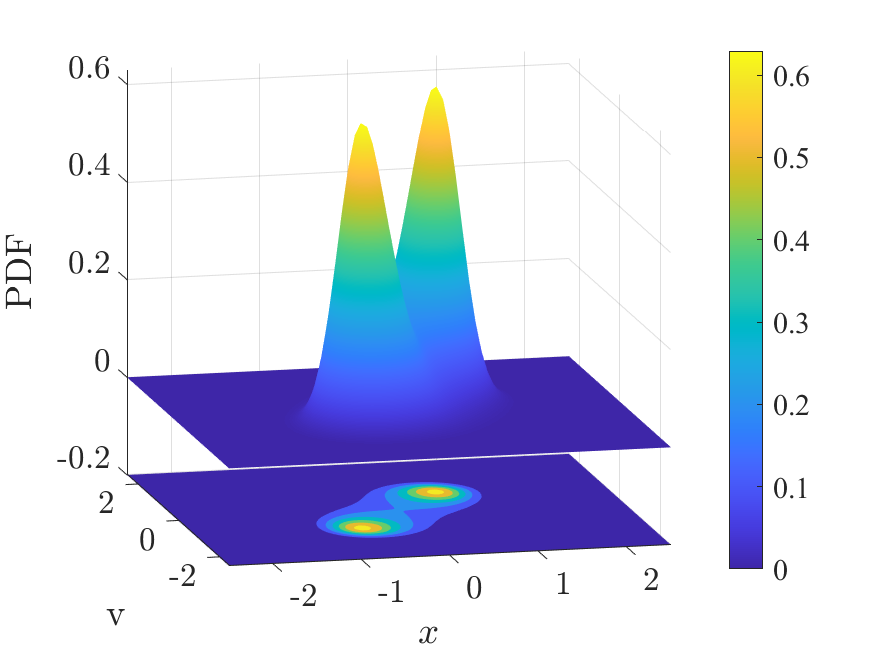}}
	\subfloat[memFPK $t=2.0$ s]{
		\includegraphics[scale=0.35]{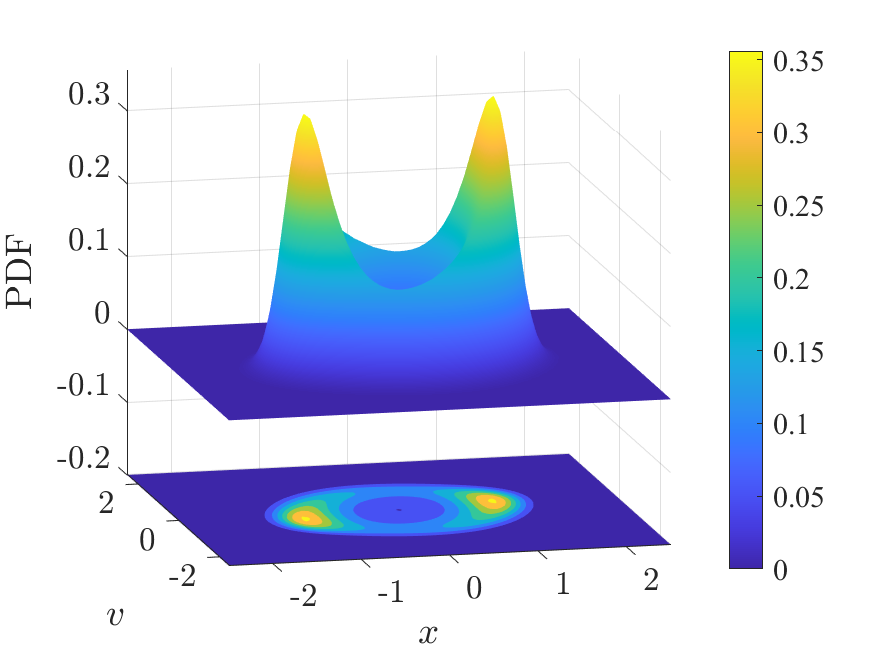}}
	\subfloat[memFPK $t=8.0$ s]{
		\includegraphics[scale=0.35]{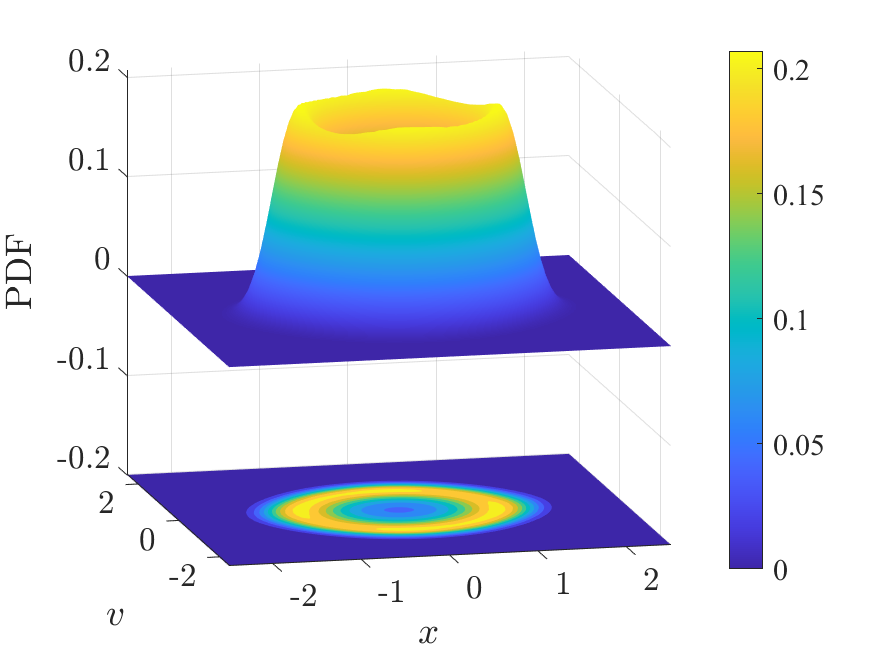}}
	\\
	\subfloat[MCS $t=1.0$ s]{
		\includegraphics[scale=0.35]{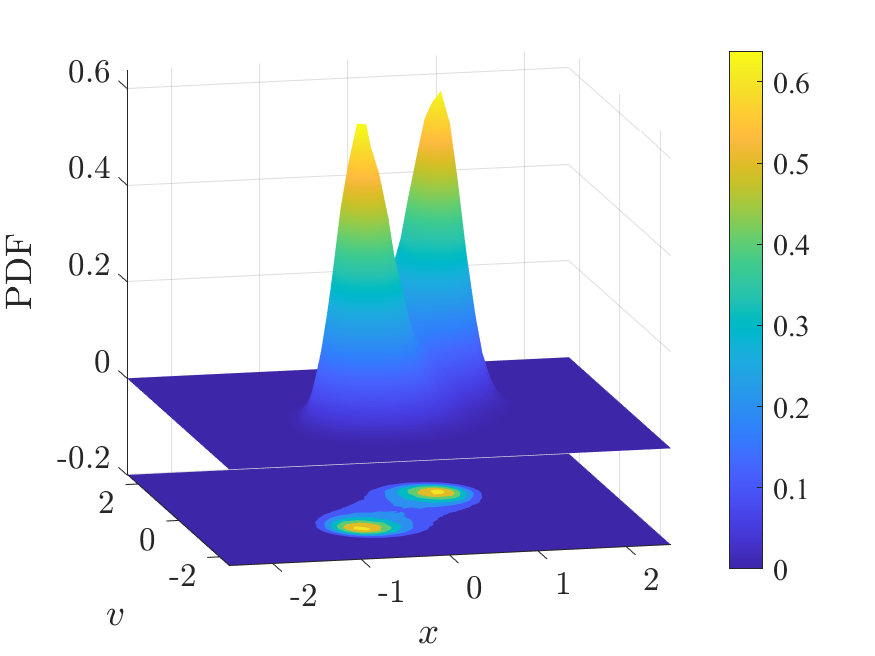}}
	\subfloat[MCS $t=2.0$ s]{
		\includegraphics[scale=0.35]{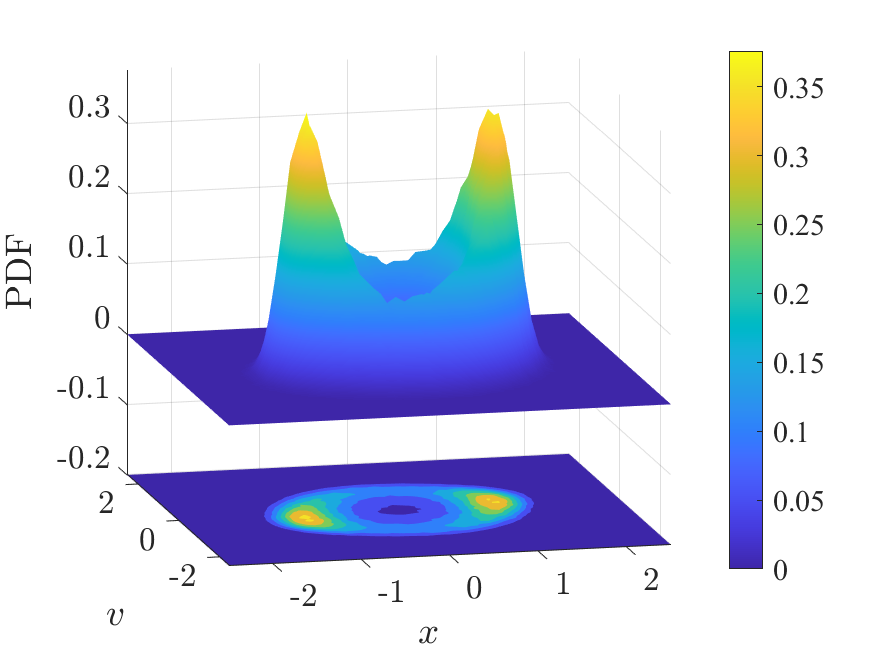}}
	\subfloat[MCS $t=8.0$ s]{
		\includegraphics[scale=0.35]{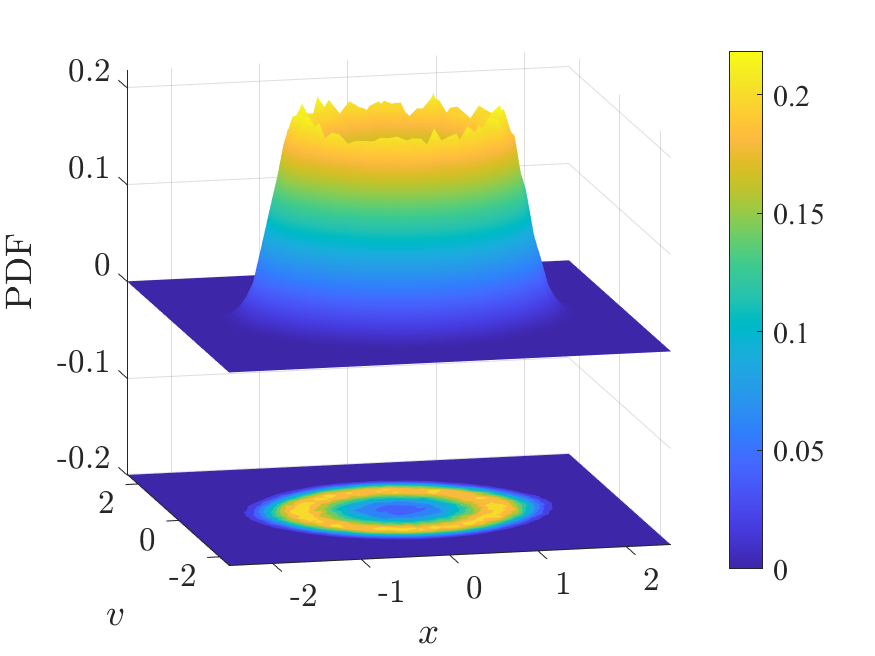}}
	\caption{Joint PDF of Example 3 obtained from memFPK and MCS solutions ($ 6\times10^6 $ samples): (a)-(c) memFPK equation solutions vs (d)-(f) MCS solutions.}
	\label{fig3-1}
\end{figure}

\begin{figure}[t!]
	\centering
	\subfloat[$X(t)$ linear scale]{
		\includegraphics[scale=0.5]{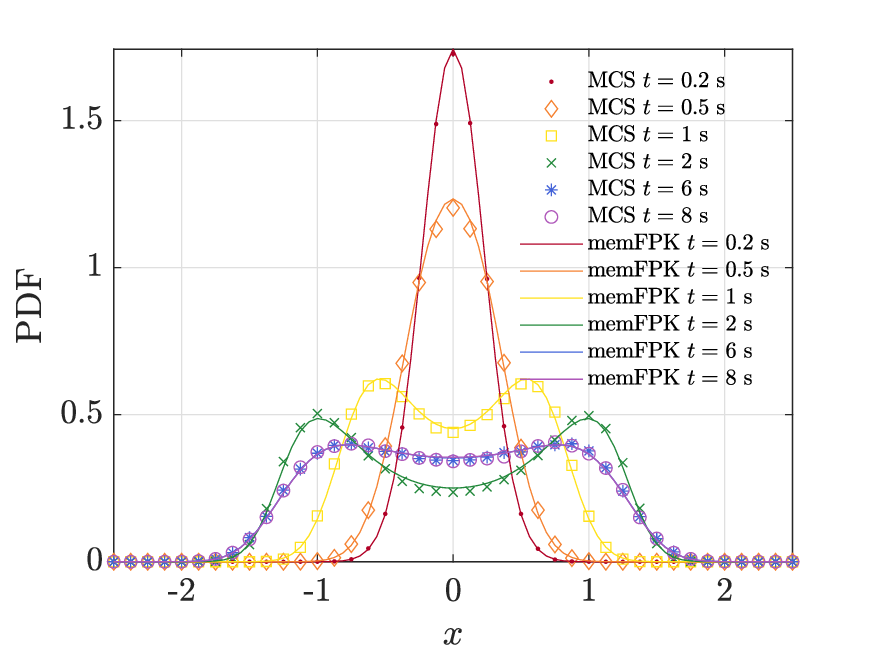}}
	\subfloat[$X(t)$ logarithmic scale]{
		\includegraphics[scale=0.5]{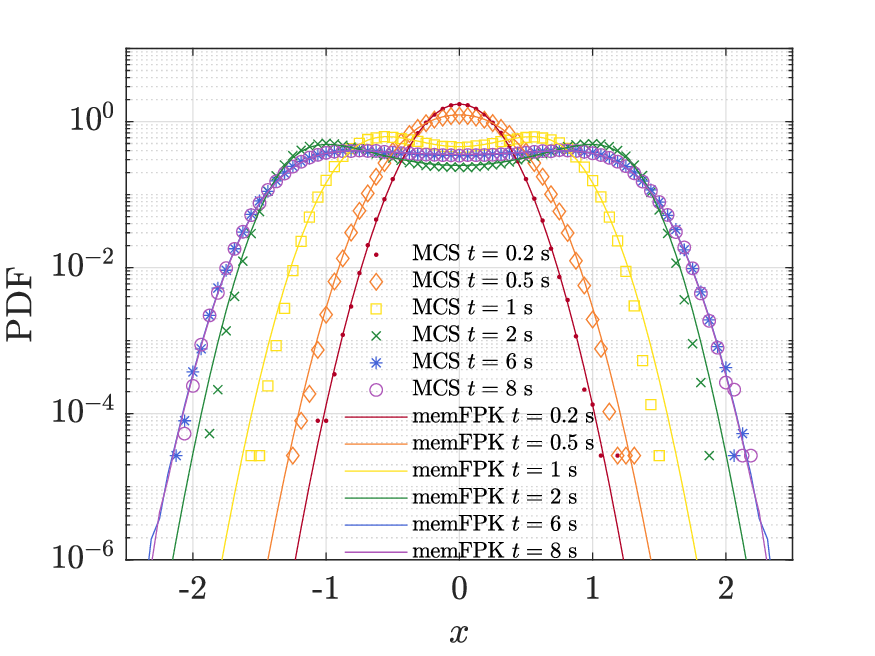}}
	\\
	\subfloat[$V(t)$ linear scale]{
		\includegraphics[scale=0.5]{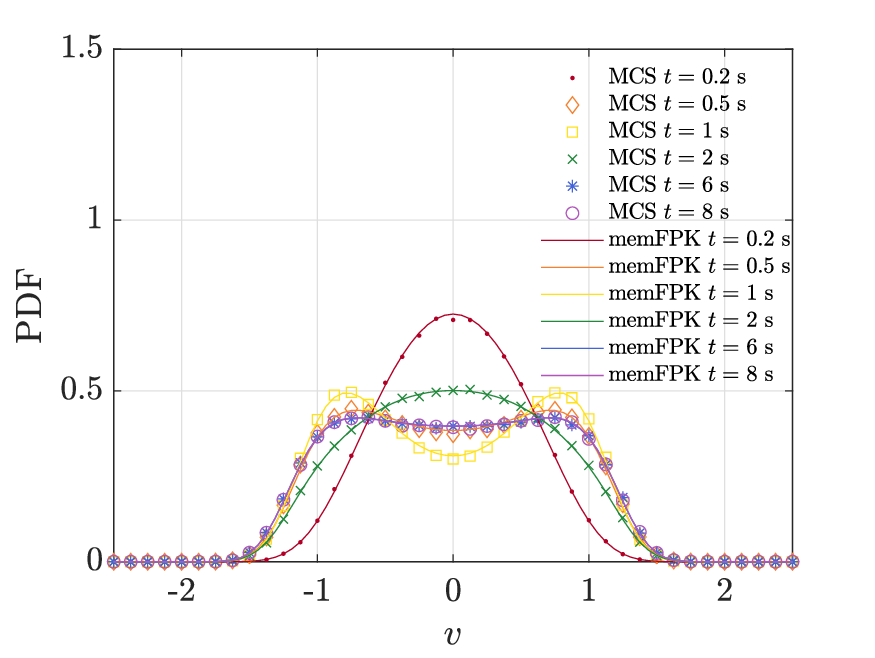}}
	\subfloat[$V(t)$ logarithmic scale]{
		\includegraphics[scale=0.5]{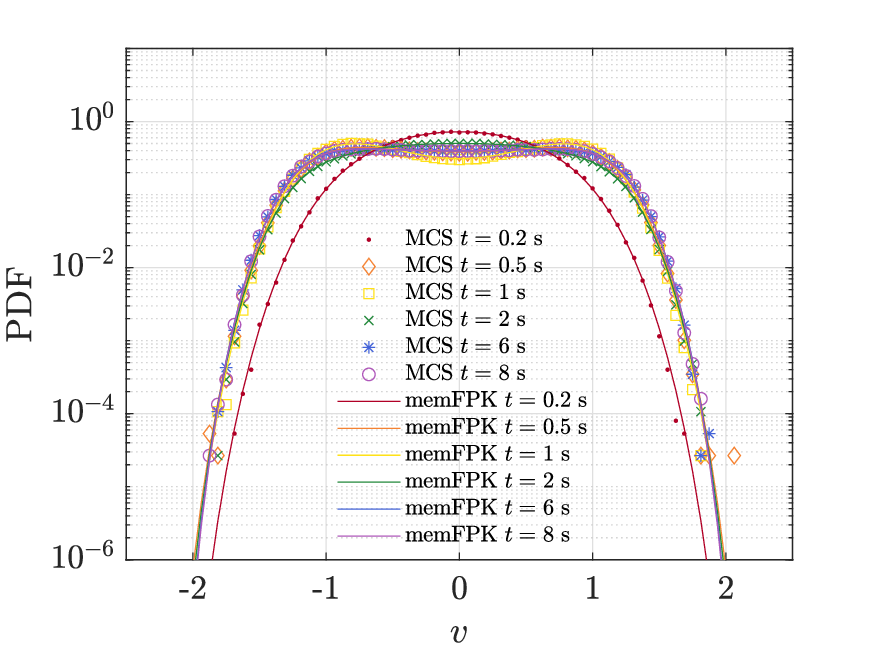}}
	\caption{Transient marginal PDFs of $X(t)$ and $V(t)$ of Example 3 obtained from memFPK and MCS solutions ($ 6\times10^6 $ samples). }
	\label{fig3-3}
\end{figure}

\begin{figure}[t!]
	\centering
	\subfloat[$X(t)$]{
		\includegraphics[scale=0.45]{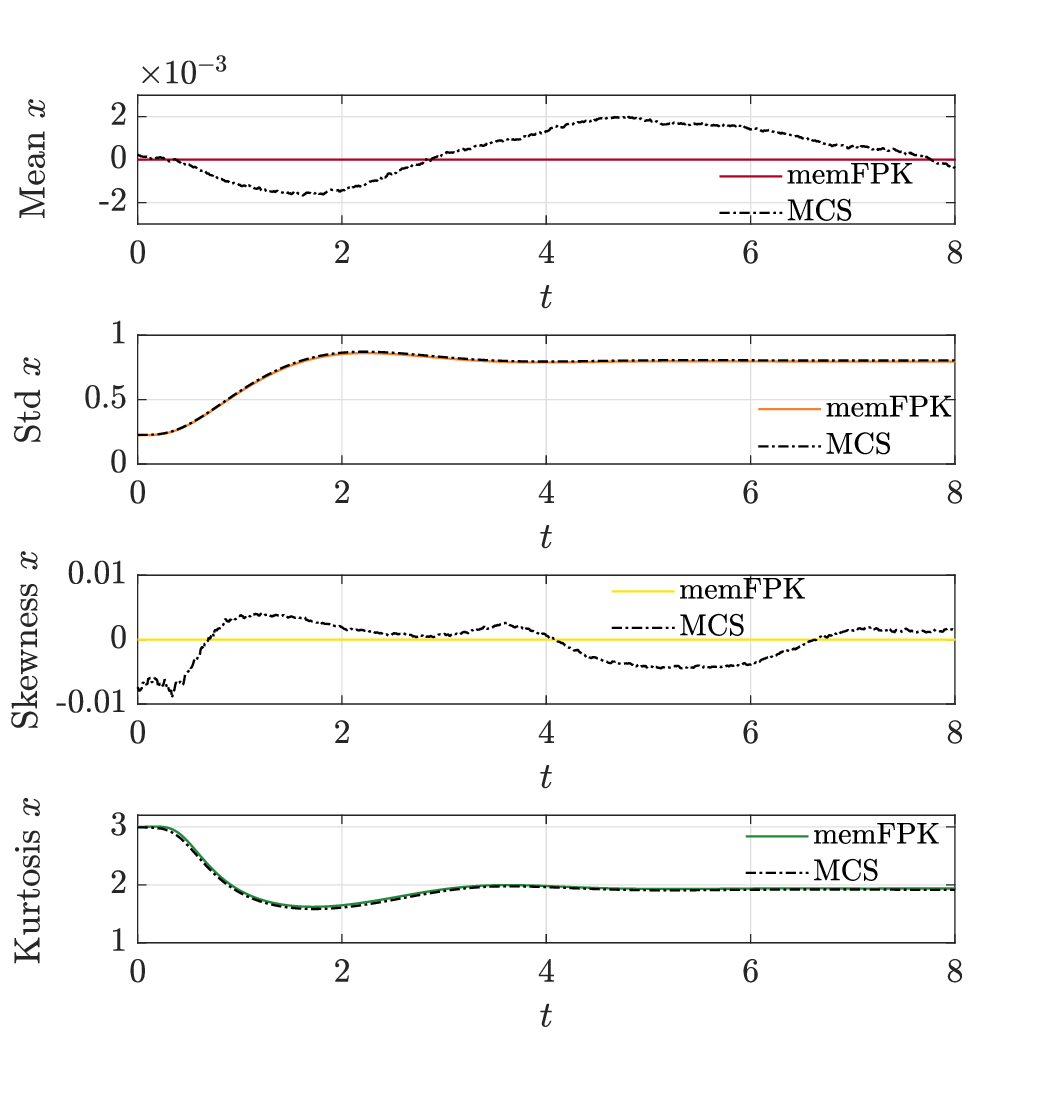}}
	\subfloat[$V(t)$]{
		\includegraphics[scale=0.45]{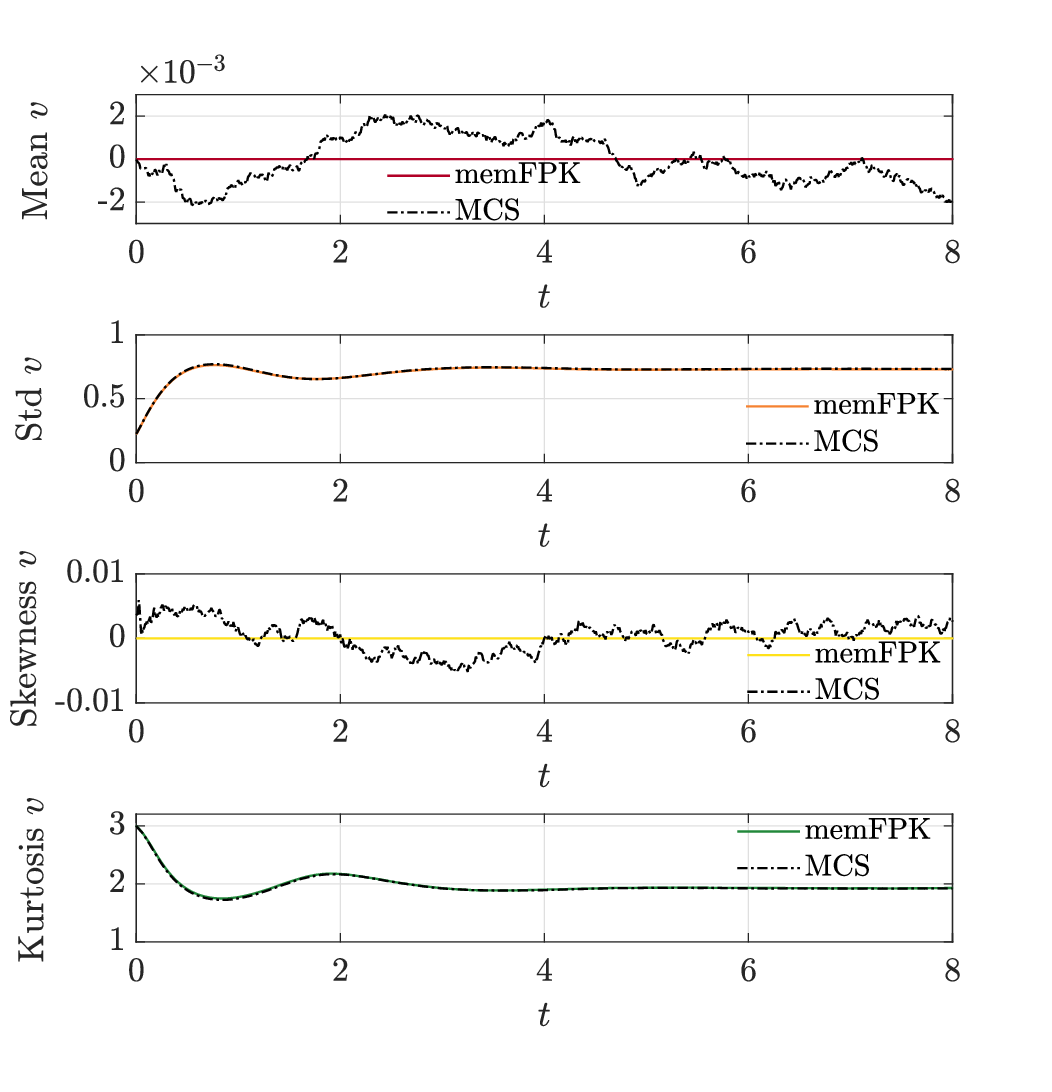}}
	\caption{Temporal evolution of the statistical moments of $X(t)$ and $V(t)$ of Example 3 obtained from memFPK and MCS solutions ($ 6\times10^6 $ samples). }
	\label{fig3-5}
\end{figure}

The resulting from solutions of the memFPK equation are subsequently compared with the results obtained from $6\times10^6$ MCS samples. In terms of PDF level, Fig.~\ref{fig3-1} shows the comparison of the joint PDF surfaces and corresponding contour plots by both the memFPK and MCS methods at representative time instants $t =$ 1.0, 2.0, and 8.0 s. The comparisons of the transient marginal PDF of $X(t)$ and $V(t)$ under both linear and logarithmic coordinate scales are shown in Fig.~\ref{fig3-3}. Regarding statistical moments level, the first four statistical moments of $X(t)$ and $V(t)$, including mean, standard deviation, skewness, and kurtosis, obtained from two methods are illustrated in Fig.~\ref{fig3-5}.

As shown in Figs.~\ref{fig3-1}-\ref{fig3-5}, the memFPK equation method effectively captures the nonlinear probabilistic characteristics of the displacement $X(t)$ and velocity $V(t)$ of the Van der Pol oscillator, and the obtained results are in good agreement with the statistics computed from $6\times10^6$ MCS samples. Unlike the Duffing oscillator, the Van der Pol oscillator is characterized by nonlinear damping, which supplies energy when the phase-space small amplitude is small and dissipates energy when it becomes large. As a result, the response probability gradually evolves from an initially concentrated distribution toward a finite amplitude oscillatory state associated with the limit cycle behavior. This behavior is clearly reflected in Fig.~\ref{fig3-1}, where the joint PDF gradually develops into a ring-like distribution in the phase plane. Accordingly, as shown in Fig.~\ref{fig3-3}, the marginal PDFs of both $X(t)$ and $V(t)$ evolve from unimodal to bimodal forms. These bimodal marginal PDFs can be understood as the projections of the ring-like phase-space distribution onto the displacement and velocity axes. Furthermore, Figs.~\ref{fig3-3}(b) and (d) show that the proposed memFPK equation method maintains excellent agreement with the MCS results even in the low-probability tail regions of the PDF values decrease to the order of $10^{-6}$. 

The temporal evolution of the first four statistical moments shown in Fig.~\ref{fig3-5}  provide further support for this probabilistic behavior. The mean values and skewness of both $X(t)$ and $V(t)$ remain close to zero, indicating that the response distributions remain symmetric. Meanwhile, the standard deviations first increase and then approach bounded levels, while the kurtosis decreases from the initial Gaussian value and remains below it during the later evolution. These trends reflect the limit cycle oscillations, further confirming the emergence of non-Gaussian response characteristics.

\subsection{Example 4}
We consider two-dimensional genetic toggle switch model proposed by Gardner \cite{gardner2000construction}
\begin{align}\label{Ex4}
	&\dot{Y_1}=\dfrac{\alpha_1}{1+Y_2^{m_1}}-Y_1+\sigma_1\xi_1^{H_1}(t),\cr
	&\dot{Y_2}=\dfrac{\alpha_2}{1+Y_1^{m_2}}-Y_2+\sigma_2\xi_2^{H_2}(t),
\end{align}
where $y_1$ and $y_2$ represent the concentrations of LacI protein and $\lambda$CI protein respectively; $\alpha_1$ and $\alpha_2$ denote the dimensionless synthesis rates in the absence of inhibitory factors; $m_1$ and $m_2$ are Hill coefficients; $\xi_1^{H_1}(t),\xi_2^{H_2}(t) $ are two independent unit FGNs, with the noise intensities $\sigma_1,\sigma_2$, respectively. The initial condition $\boldsymbol{Y}(0)=\boldsymbol{y}_0=(y_{10},y_{20})^T$ is random vector that follows the joint Gaussian distribution  $N(\boldsymbol{\mu}_{0},\Sigma_{0})$.

Set the state vector $\boldsymbol{Y}=(Y_1,Y_2)^{T}$, Eq.~(\ref{Ex4}) can be expressed in vector form, where
\begin{align*}
	\boldsymbol{\hat f}(\boldsymbol{Y})=\begin{pmatrix}\dfrac{\alpha_1}{1+Y_2^{m_1}}-Y_1 \\ \dfrac{\alpha_2}{1+Y_1^{m_2}}-Y_2\end{pmatrix},
	\hat \Sigma=\begin{pmatrix}  \sigma_1&0\\0&\sigma_2 \end{pmatrix},
	\boldsymbol{\xi}^{\boldsymbol{H}}=\begin{pmatrix}  \xi_1^{H_1} \\ \xi_2^{H_2} \end{pmatrix}.
\end{align*}
Then, refer to Theorem \ref{them-1}, the response joint PDF of system~(\ref{Ex4}) satisfies the memFPK equation \eqref{memFPK-2d-nonlinear-1} with the memory-dependent drift and diffusion coefficients (\ref{2d-coffs-1}), where
$$\nabla\boldsymbol{\hat f}(\boldsymbol{Y})=\begin{pmatrix}
	-1 & \dfrac{\alpha_1 m_1 Y_2^{m_1-1}}{(1+Y_2^{m_1})^2} \\
	\dfrac{\alpha_2 m_2 Y_1^{m_2-1}}{(1+Y_1^{m_2})^2} & -1
\end{pmatrix}.$$ 

In this case, the DLM treatment estimates memory-dependent diffusion coefficients using 2000 MCS samples. The results are substituted in the memFPK equation, which is solved numerically with FD method. The computational domains are set at $y_1 \in [-2, 5]$ and $y_2 \in [-2, 5]$, with spatial steps $\Delta y_1 = 0.1$ and $\Delta y_2 = 0.1$, and time step $\Delta t=$ 0.001 s. And the parameters used are $\alpha_1=\alpha_2= 2.5, m_1=m_2=2, \sigma_1^2 = 0.25, \sigma_2^2 = 0.36, \boldsymbol{\mu}_{0}=(1.2,1)^T,\Sigma_{0}={\rm diag}(0.05,0.05), H_1=0.8, H_2=0.7$.

\begin{figure}[t!]
	\centering
	\subfloat[memFPK $t=$ 2.0 s]{
		\includegraphics[scale=0.35]{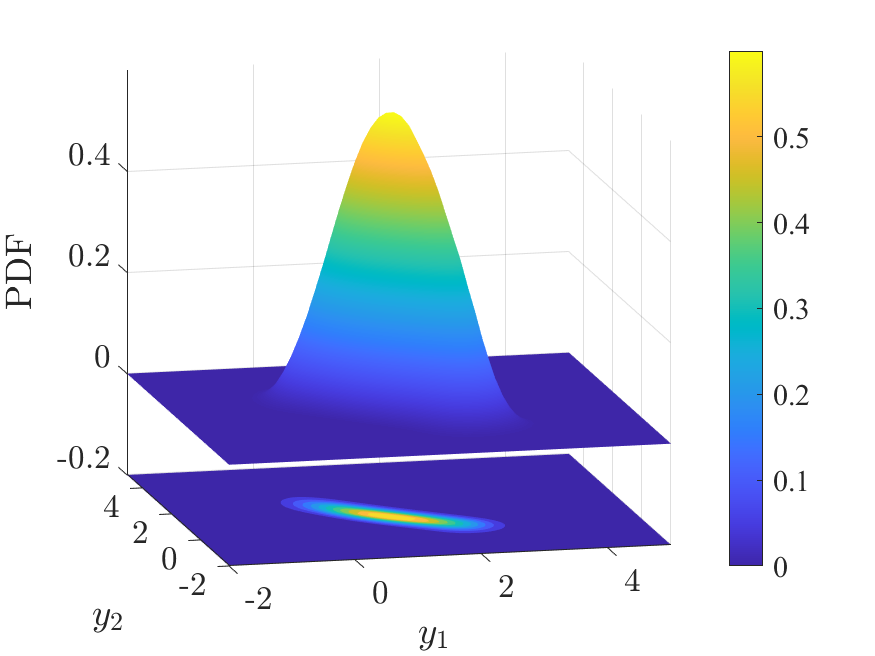}}
	\subfloat[memFPK $t=$ 5.0 s]{
		\includegraphics[scale=0.35]{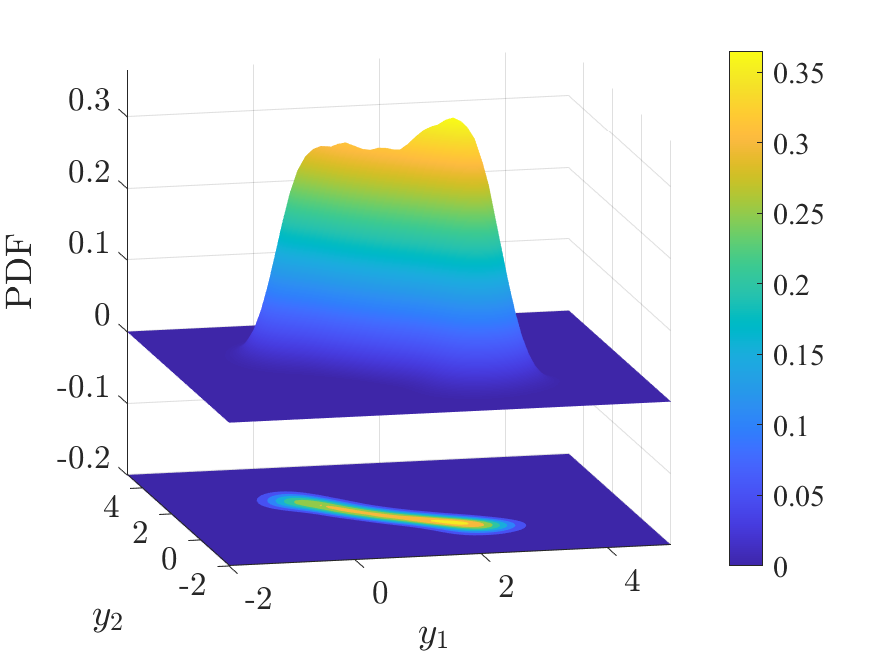}}
	\subfloat[memFPK $t=$ 16.0 s]{
		\includegraphics[scale=0.35]{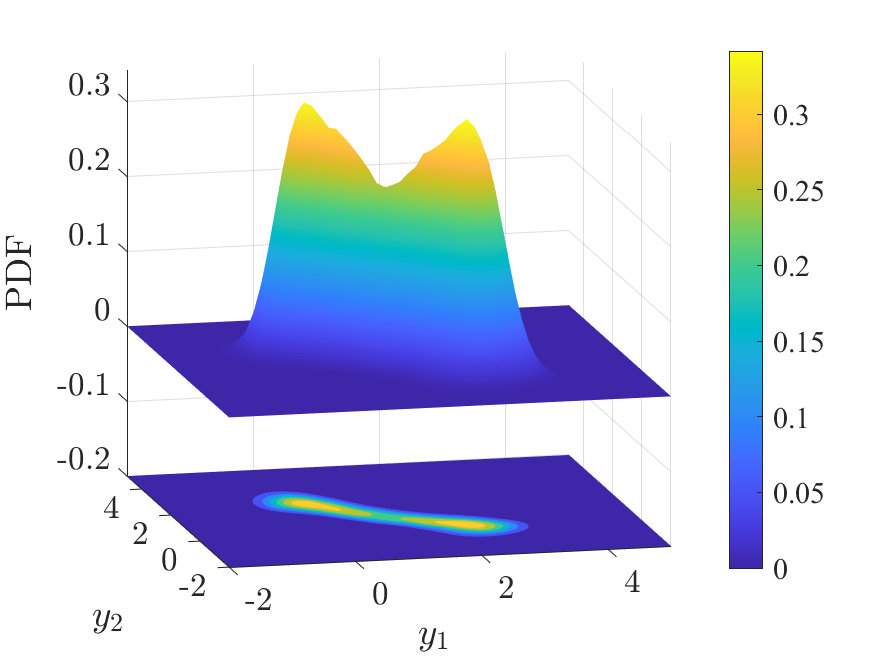}}
	\\
	\subfloat[MCS $t=$ 2.0 s]{
		\includegraphics[scale=0.35]{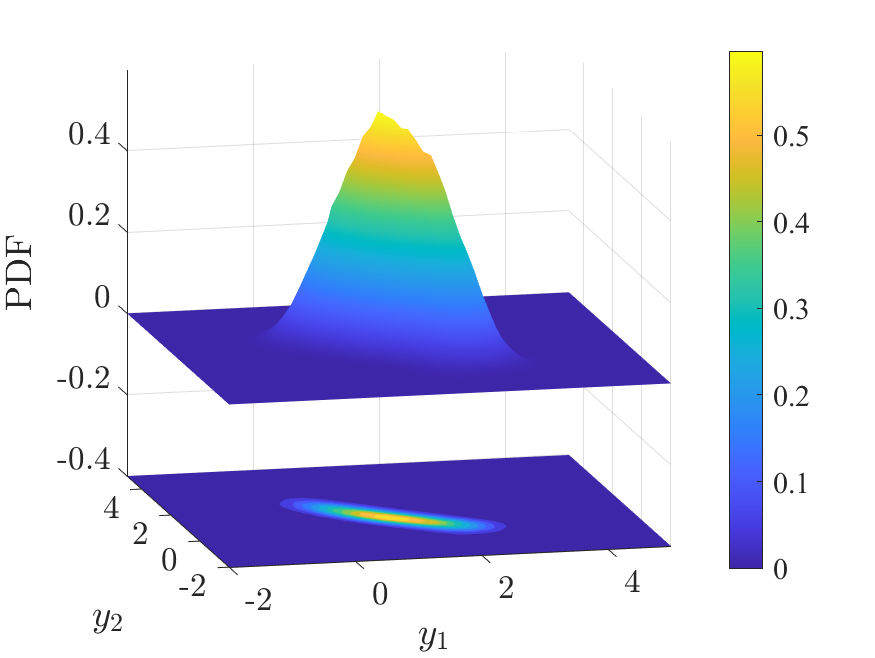}}
	\subfloat[MCS $t=$ 5.0 s]{
		\includegraphics[scale=0.35]{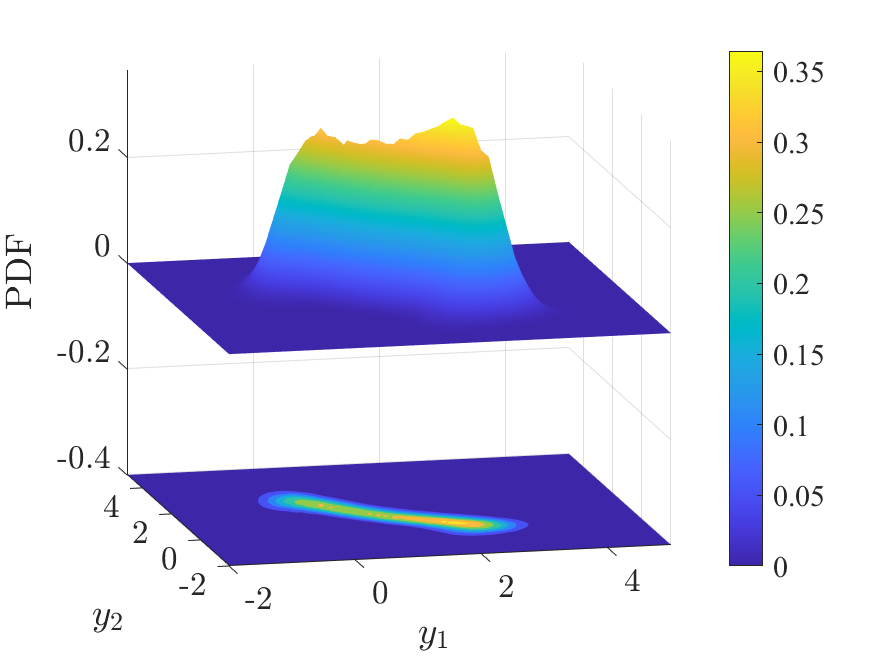}}
	\subfloat[MCS $t=$ 16.0 s]{
		\includegraphics[scale=0.35]{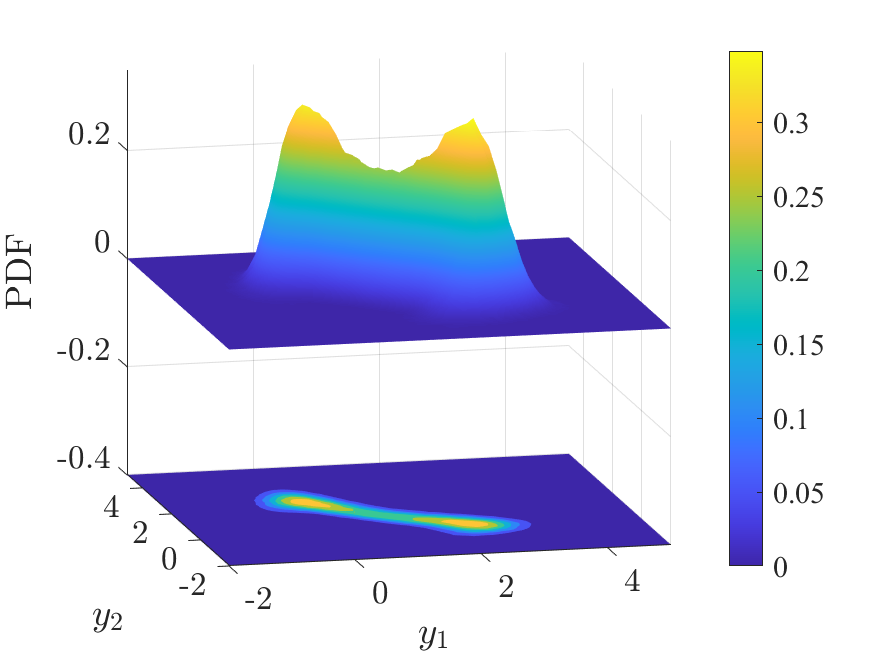}}
	\caption{Joint PDF of Example 4 obtained from memFPK and MCS solutions ($ 6\times10^6 $ samples): (a)-(c) memFPK equation solutions vs (d)-(f) MCS solutions.}
	\label{fig4-1}
\end{figure}

\begin{figure}[t!]
	\centering
	\subfloat[$Y_1(t)$ linear scale]{
		\includegraphics[scale=0.5]{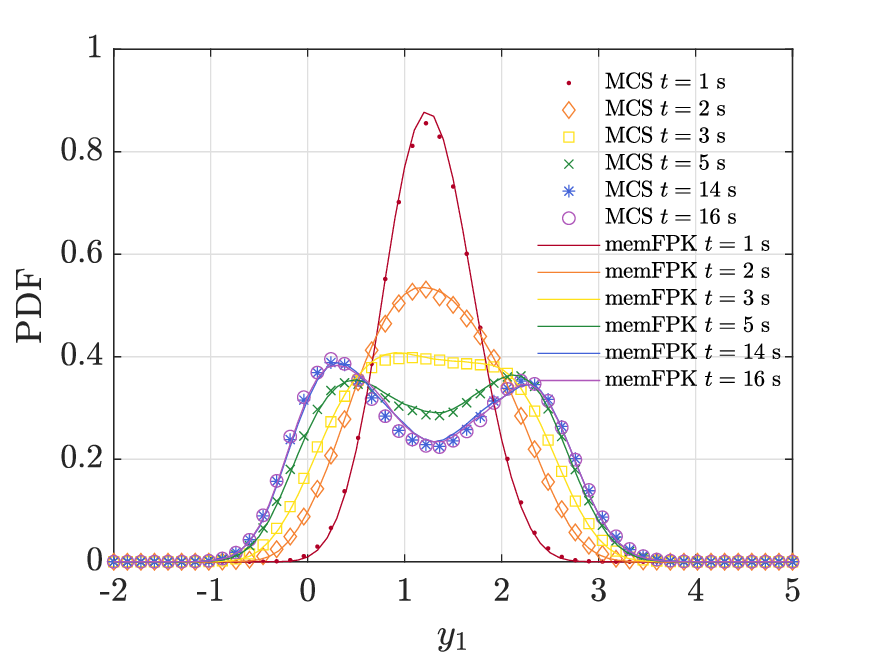}}
	\subfloat[$Y_1(t)$ logarithmic scale]{
		\includegraphics[scale=0.5]{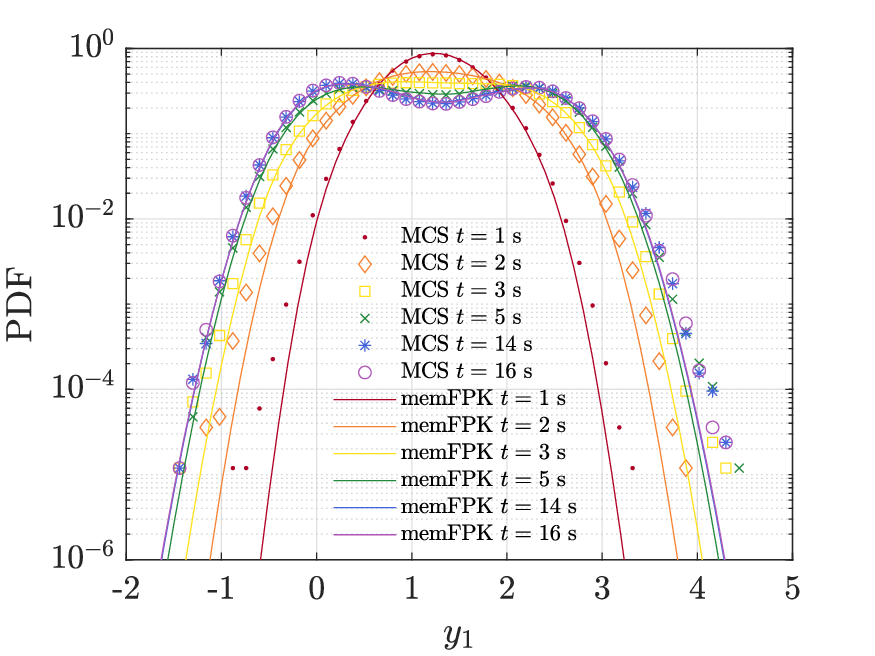}}
	\\
	\subfloat[$Y_2(t)$ linear scale]{
		\includegraphics[scale=0.5]{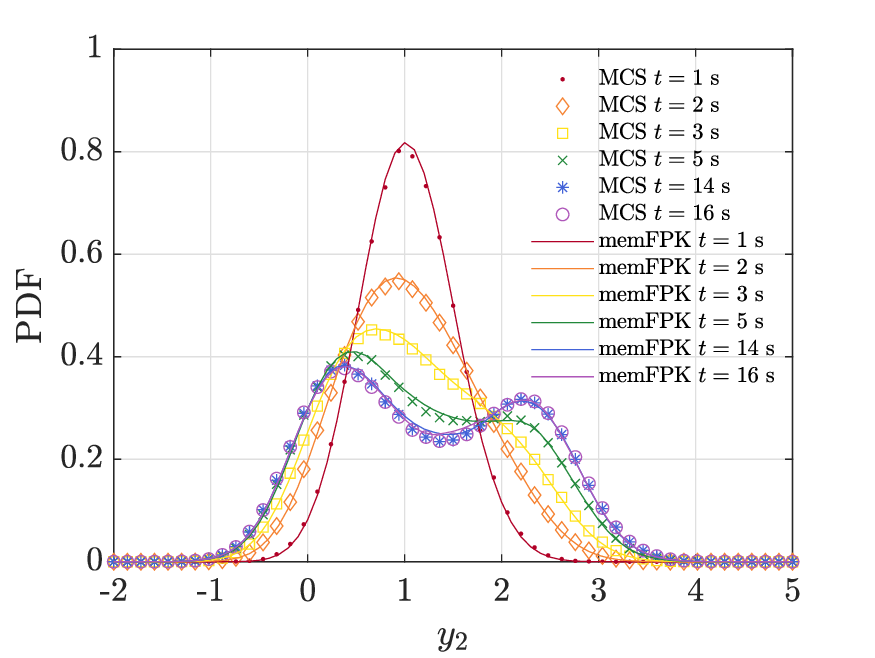}}
	\subfloat[$Y_2(t)$ logarithmic scale]{
		\includegraphics[scale=0.5]{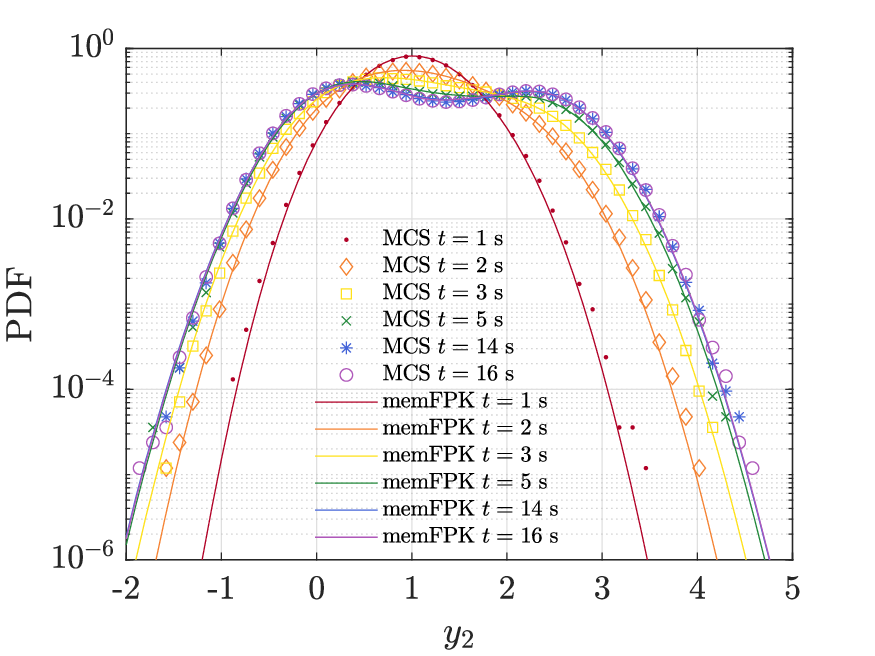}}
	\caption{Transient marginal PDFs of $Y_1(t)$ and $Y_2(t)$ of Example 4 obtained from memFPK equation and MCS solutions ($ 6\times10^6 $ samples). }
	\label{fig4-3}
\end{figure}

\begin{figure}[t!]
	\centering
	\subfloat[$Y_1(t)$]{
		\includegraphics[scale=0.45]{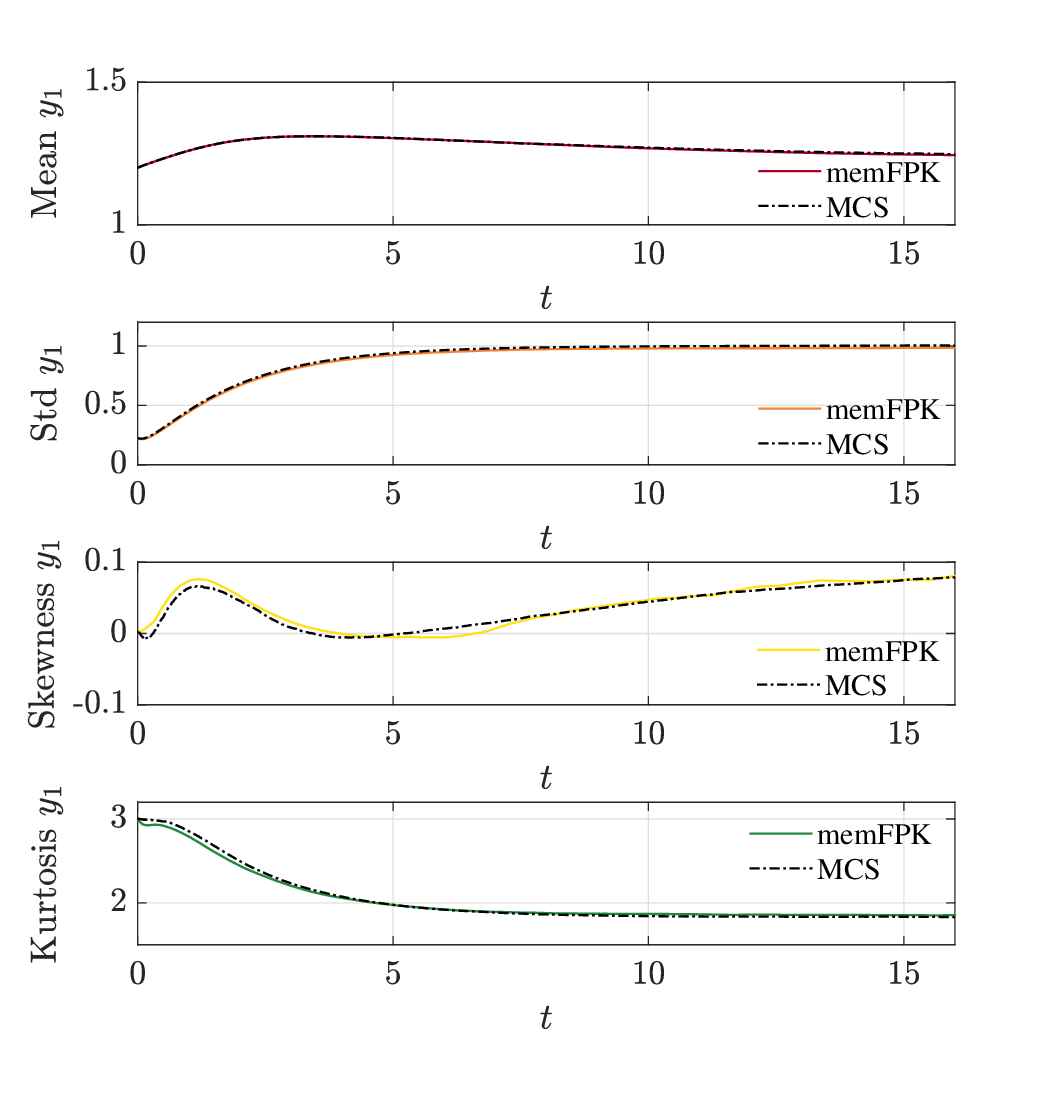}}
	\subfloat[$Y_2(t)$]{
		\includegraphics[scale=0.45]{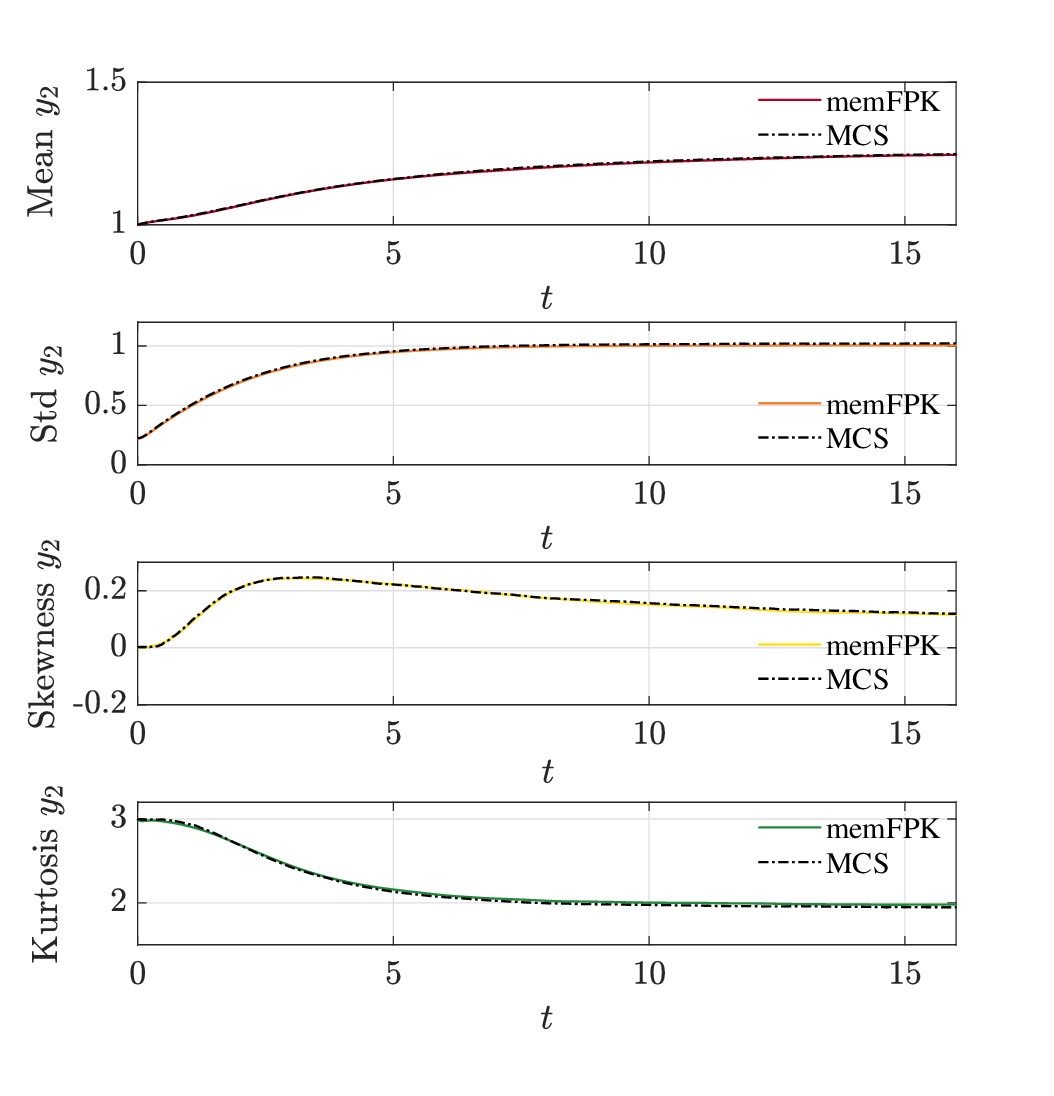}}
	\caption{Temporal evolution of the statistical moments of $Y_1(t)$ and $Y_2(t)$ of Example 4 obtained from memFPK equation and MCS solutions ($ 6\times10^6 $ samples).}
	\label{fig4-5}
\end{figure}

The numerical results were validated by comparing with $6\times10^6$ MCS samples. Fig.~\ref{fig4-1} presents the joint PDF at representative time instants $t =$ 2.0, 5.0, and 16.0 s, using both surface and contour plots. The corresponding marginal PDF of $Y_1(t)$ and $Y_2(t)$, plotted on linear and logarithmic scales, are presented in Fig.~\ref{fig4-3}, while the time evolutions of the first four statistical moments are illustrated in Fig.~\ref{fig4-5}. As shown in these figures, the memFPK equation solutions are in good agreement with the statistics obtained from $6\times10^6$ MCS samples.

It is notable that, under the selected parameter settings, Figs.~\ref{fig4-1} and \ref{fig4-3} exhibit pronounced asymmetric bimodal structure in both $Y_1(t)$ and $Y_2(t)$. These features are consistent with the mutual repression mechanism of the genetic toggle switch. Under FGN excitation, the observed bimodal responses suggest stochastic switching behavior between these states. The statistical moment evolutions in Fig.~\ref{fig4-5} further characterize the transient probabilistic behavior of the genetic toggle switch. The mean values of $Y_1(t)$ and $Y_2(t)$ evolve toward stable levels, while their standard deviations increase and then approach bounded values, indicating the gradual spreading and stabilization of the response distributions. In addition, the positive skewness indicates the asymmetric response distributions, and the kurtosis gradually decreases below the Gaussian value 3, which is consistent with the formation of bimodal non-Gaussian probability structures. These results demonstrate that the memFPK equation method can accurately capture the asymmetric switching behavior and non-Gaussian characteristics of the genetic regulatory system.

\section{Discussions and conclusions}\label{sec13}
This paper presents a framework for analyzing the transient probabilistic responses of two-dimensional nonlinear stochastic dynamical systems under external FGN excitation, with a special focus on nonlinear SDOF oscillators. By using the FWIS integral theory, we establish a new memFPK equation to describe the transient evolution of the joint PDF of two-dimensional systems excited by FGN. To solve this equation numerically, we develop a data-driven DLM treatment to estimate the memory-dependent diffusion coefficients that involve conditional expectations. After estimating these coefficients, we integrate them into the memFPK equation, which is then solved using a FD scheme to obtain the transient probabilistic responses. We verify the validity and effectiveness of the proposed framework through numerical examples. The results show that the method can accurately capture transient joint PDFs, marginal PDFs, low-probability tail behavior, and statistical moments. Comparisons with analytical solutions or MCS results further prove the accuracy and numerical reliability of the proposed approach. These findings confirm that the method is capable of describing the transient probabilistic responses induced by long-memory stochastic excitation.

This study offers an effective method for the probabilistic analysis of non-Markovian stochastic dynamics in engineering and biomedical fields. In practical engineering research, the effectiveness of the proposed memFPK framework has been verified in the probabilistic analysis of tilt-rotor aircraft whirl flutter \cite{kong2026analytical}. In biomedical modeling, the numerical example of genetic toggle switches presented in this work proves that the established method is capable of characterizing non-Gaussian features and stochastic switching phenomena in nonlinear regulatory systems. In a general sense, this method broadens the application scope of traditional FPK equations beyond classic Markovian systems, and provides a straightforward approach to investigate transient joint PDFs, marginal PDFs, tail characteristics and various non-Gaussian response behaviors under FGN excitation.

For future research, we will extend the proposed memFPK equation method to multi-degree-of-freedom nonlinear systems under simultaneous external and parametric FGN excitations. Such extensions will involve more demanding calculations and may further increase the difficulty of estimating memory-dependent coefficients in sparsely sampled regions. To address these challenges, more efficient numerical algorithms will be explored for high-dimensional non-Markovian probabilistic analysis. Recent studies \cite{zhang2022solving,zhang2023deep} have shown that deep learning can provide useful tools for solving FPK equations, while Wang et al. \cite{wang2025pseudo} further extended this idea to parameterized FPK equations and achieved notable improvements in computational efficiency. These developments suggest that deep learning-based methods may also help improve both the numerical solution of high-dimensional memFPK equations and the estimation of memory-dependent coefficients under limited sample information. These above efforts can further broaden the practical application of this theory in engineering analysis and structural design.

\section*{CRediT authorship contribution statement}
\textbf{Lifang Feng:} Conceptualization, Formal analysis, Investigation, Methodology, Software, Visualization, Validation, Writing-original draft, Writing-review \& editing.
\textbf{Bin Pei:} Conceptualization, Funding acquisition, Methodology, Writing-original draft, Writing-review \& editing.
\textbf{Yong Xu:} Conceptualization, Funding acquisition, Project administration, Supervision, Writing-review \& editing.

\section*{Declaration of competing interest}
The authors declare that they have no known competing financial interests or personal relationships that could have appeared to influence the work reported in this paper.

\section*{Acknowledgements}
Y. Xu was supported by the National Natural Science Foundation (NSF) of China under Grant No. U2441204. B. Pei and L. Feng  were supported by NSF of China under Grant No. 12172285 and Fundamental Research Funds for the Central Universities.

\appendix
\renewcommand{\thetheorem}{\Alph{section}.\arabic{theorem}}

\setcounter{theorem}{0}

\section{Mathematical concepts and definitions}\label{app-sec-1}
Several basic concepts required for the derivation of the memFPK equation are briefly introduced, including fractional Brownian motion (FBM), FGN, the symmetric pathwise integral, the FWIS integral, and the fractional It\^o's formula.

\subsection{FBM and FGN}
For a Hurst index $H\in(0,1)$, a one-dimensional FBM $B^{H}(t)$ is a centered Gaussian process satisfying the following properties
\begin{align*}
	\mathbb{E}[B^{H}(t)] = 0, \qquad
	\mathbb{E}[B^{H}(t)B^{H}(s)]=\frac{1}{2}(t^{2H}+s^{2H}-|t-s|^{2H}),
\end{align*}
for $s,t\geq 0$. When $H=1/2$, $B^{H}(t)$ reduces to standard Brownian motion. The increment $B^{H}(t)-B^{H}(s)$ is Gaussian distributed with zero mean and variance
\begin{align*}
	\mathbb{E}[B^{H}(t)-B^{H}(s)]= 0, \qquad
	\mathbb{E}[(B^{H}(t)-B^{H}(s))^2]=\vert t-s\vert^{2H},\qquad s,t\geq 0.
\end{align*}
For $H>1/2$, the increments are positively correlated, which gives rise to long-range dependence.

Similar to the relationship between standard Brownian motion and GWN, a unit FGN $\xi^{H}(t)$ is formally defined as the generalized time derivative of FBM, i.e.,
\begin{align}\label{eq-app-3}
	\xi^{H}(t)=\frac{\mathrm{d}B^{H}(t)}{\mathrm{d}t},\quad {\rm or}\quad B^{H}(t)=\int_{0}^{t}\xi^{H}(s)\mathrm{d}s.
\end{align}

\subsection{Stochastic integrals for FBM}
For a one-dimensional nonlinear stochastic system driven by FGN,
\begin{align}\label{eq-app-6}
	\dot{Y}(t)=f(Y(t))+g(Y(t))\xi^{H}(t).
\end{align}
Using the relation between FBM and FGN in Eq.~\eqref{eq-app-3}, Eq.~\eqref{eq-app-6} can be rewritten as the following SDE in the sense of symmetric pathwise integral (Stratonovich-type SDE)
\begin{align}\label{eq-app-7}
	\mathrm{d}Y(t)=f(Y(t))\mathrm{d}t+
	g(Y(t))\mathrm{d}^{\circ}B^{H}(t).
\end{align}
Here, $\mathrm{d}^{\circ}$ denotes the so-called symmetric pathwise integral \cite{nualart2003stochastic}, which is defined as
\begin{align}\label{eq-app-8}
	\int_0^t g(s)\mathrm{d}^{\circ}B^{H}(s)
	=\lim_{|\Pi|\to 0}\sum_{i=0}^{n-1}\frac{g(t_i)+g(t_{i+1})}{2}(B^{H}(t_{i+1})-B^{H}(t_i)),
\end{align}
where $\Pi=\{0=t_0<t_1<\cdots<t_n=t\}$ is a partition of $[0,t]$. For $H=1/2$, Eq.~\eqref{eq-app-8} reduces to the Stratonovich integral, so the symmetric pathwise integral is also called the Stratonovich-type integral.

Besides the symmetric pathwise integral, forward and backward pathwise integrals with respect to FBM can also be defined \cite{biagini2008stochastic,zahle2002forward}. Under suitable regularity conditions, these pathwise integrals coincide for $H>1/2$ \cite{lu2017stationary}, and thus the pathwise formulation provides a natural interpretation of Eq.~\eqref{eq-app-7}. However, such pathwise integrals do not, in general, possess the zero-mean property required when taking expectations in the derivation of the FPK equation from the It\^o's formula. Therefore, it is more convenient to work in the FWIS framework for the subsequent analysis.

The FWIS integral is a divergence-type stochastic integral with respect to the FBM for Hurst index $1/2<H<1$ by using the fractional white-noise analysis method \cite{biagini2008stochastic}. For an admissible integrand $g(t)$, it can be written as
\begin{align*}
	\int_0^t g(s)\mathrm{d}^{\diamond}B^{H}(s)
	=
	\lim_{|\Pi|\to 0}\sum_{i=0}^{n-1}g(t_i)\diamond(B^{H}(t_{i+1})-B^{H}(t_i)),
\end{align*}
where $\diamond$ denotes the Wick product and $\mathrm{d}^{\diamond}$ denotes the FWIS integral \cite[Definition 3.4.1.]{biagini2008stochastic}. A key property of this integral is that \cite[Theorem 3.6.1.]{biagini2008stochastic}, $
	\mathbb{E}\Big[\int_0^t g(s)\mathrm{d}^{\diamond}B^{H}(s)\Big]=0, $
whenever the integral is well defined, so the FWIS integral is also called the It\^o-type integral.

For sufficiently regular integrands, the symmetric pathwise and FWIS integrals are related by \cite[Theorem 3.12.]{ducan2000stochastic}
\begin{align*}
	\int_0^t g(s)\mathrm{d}^{\circ}B^{H}(s)
	=
	\int_0^t g(s)\mathrm{d}^{\diamond}B^{H}(s)
	+
	\int_0^t D_s^{H}(g(s))\mathrm{d}s,
\end{align*}
where $ D_s^{H}(g(s))$ is called the Malliavin derivative of $g(s)$ associated with the FBM $B^{H}$. Further details on the Malliavin derivative can be found in \cite{ducan2000stochastic,hu2003fractional}.

In the present two-dimensional setting, we denote
\begin{align*}
	\boldsymbol{B}^{\boldsymbol{H}}(t)
	=\big(B_1^{H_1}(t),B_2^{H_2}(t)\big)^T,
	\qquad
	\boldsymbol{\xi}^{\boldsymbol{H}}(t)
	=\big(\xi_1^{H_1}(t),\xi_2^{H_2}(t)\big)^T,
\end{align*}
where $B_1^{H_1}(t)$ and $B_2^{H_2}(t)$ are independent FBMs with $H_i\in(1/2,1)$, $i=1,2$, and
\begin{align*}
	\frac{\mathrm{d}\boldsymbol{B}^{\boldsymbol{H}}(t)}{\mathrm{d}t}
	=\boldsymbol{\xi}^{\boldsymbol{H}}(t), \qquad
	\mathbb{E}\big[\boldsymbol{B}^{\boldsymbol{H}}(t)\big]
	=\boldsymbol{0}, \qquad
	\mathbb{E}\big[\boldsymbol{B}^{\boldsymbol{H}}(t)(\boldsymbol{B}^{\boldsymbol{H}}(t))^T\big]
	={\rm diag}\big(t^{2H_1},t^{2H_2}\big).
\end{align*}

Consider the two-dimensional nonlinear stochastic system driven by FGNs,
\begin{align}\label{eq-app-61}
	\dot{\boldsymbol{Y}}(t)
	=
	\boldsymbol{\hat f}(\boldsymbol{Y}(t))
	+
	\hat{\Sigma}\boldsymbol{\xi}^{\boldsymbol{H}}(t),
\end{align}
where $$\hat{\Sigma}=
\begin{pmatrix}
	\sigma_{11} & \sigma_{12}\\
	\sigma_{21} & \sigma_{22}
\end{pmatrix},$$
is a constant noise-intensity matrix. Similar to the one-dimensional case, Eq.~\eqref{eq-app-61} can be rewritten as the following Stratonovich-type SDE
\begin{align}\label{eq-app-71}
	\mathrm{d}\boldsymbol{Y}(t)
	=
	\boldsymbol{\hat f}(\boldsymbol{Y}(t))\mathrm{d}t
	+
	\hat{\Sigma}\mathrm{d}^{\circ}\boldsymbol{B}^{\boldsymbol{H}}(t).
\end{align}
 Since the stochastic excitation is additive, and the diffusion matrix $\hat{\Sigma}$ is constant. Hence, $D_s^{\boldsymbol{H}}(\hat{\Sigma})=\boldsymbol{0}$, which implies that $\int_0^t\hat{\Sigma}\mathrm{d}^{\circ}\boldsymbol{B}^{\boldsymbol{H}}(s)
	=
	\int_0^t\hat{\Sigma}\mathrm{d}^{\diamond}\boldsymbol{B}^{\boldsymbol{H}}(s). $
Therefore, the Stratonovich-type SDE \eqref{eq-app-71} is equivalent to the It\^o-type SDE
\begin{align}\label{eq-app-13}
	\mathrm{d}\boldsymbol{Y}(t)
	=
	\boldsymbol{\hat f}(\boldsymbol{Y}(t))\mathrm{d}t
	+
	\hat{\Sigma}\mathrm{d}^{\diamond}\boldsymbol{B}^{\boldsymbol{H}}(t).
\end{align}
\subsection{Several useful Lemmas}
By employing similar arguments as in \cite[Theorem 3.12, Theorem 4.2, Theorem 4.5]{ducan2000stochastic} for one-dimensional and \cite[Section 3.11]{biagini2008stochastic} for mult-dimensional cases, the following lemmas can be derived, which are essential for obtaining the memFPK equation.

\begin{lemma}[Fractional It\^o's formula]\label{mul-itoformula}
	Let $ \boldsymbol{Y}(t)=(Y_1(t),Y_2(t))^{T} $ be the solution of Eq.~(\ref{eq-app-13}). If Hurst index $ \boldsymbol{H}=(H_1,H_2) \in (1/2,1)^2$, $\boldsymbol{\hat f} = (f_1,f_2)^T$, $\hat \Sigma$ is non-zero constant matrix. Then, suppose that $ F\in C^{2,1}(\mathbb{R}^2\times\mathbb{R}) $ with bounded second-order derivatives with respect to $\boldsymbol{Y}$ and first-order derivatives with respect to $t$. Then for $ t\in [0,T] $
	\begin{align}\label{mul-itofor}
			\mathrm{d}F(\boldsymbol{Y}(t),t)
			=&\frac{\partial F}{\partial t}(\boldsymbol{Y}(t),t)\mathrm{d}t
			+\sum_{i=1}^{2}f_i(\boldsymbol{Y}(t))\frac{\partial F}{\partial y_i}(\boldsymbol{Y}(t),t) \mathrm{d}t+\sum_{i=1}^{2}\sum_{k=1}^{2}\sigma_{ik}\frac{\partial F}{\partial y_i}(\boldsymbol{Y}(t),t)\mathrm{d}^{\diamond}B_k^{H_k}(t)\cr
			&+\sum_{i=1}^{2}\sum_{j=1}^{2}\sum_{k=1}^{2}\sigma_{ik}D^{H_k}_{t}(Y_j(t))\frac{\partial^2 F}{\partial y_i\partial y_j}(\boldsymbol{Y}(t),t)\mathrm{d}t,
	\end{align}
	where $ D^{H_k}_{t} (Y_j(t)) $ is the Malliavin derivative of $Y_j(t)$ with respect to $ B_k^{H_k}(t)$ for $j,k=1,2$.
\end{lemma}
The last term in Eq.~\eqref{mul-itofor} is the main feature distinguishing the fractional It\^o's formula from the classical It\^o's formula under GWN excitation. It carries the memory contribution induced by FGN and ultimately gives rise to the memory-dependent diffusion coefficients in the memFPK equation.

\begin{lemma}\label{malian-to-fwis}
	If $ \boldsymbol{\hat f}, \hat \Sigma $ be random function vector and non-zero constant matrix such that the assumptions of Lemma \ref{mul-itoformula} are satisfied, then for $r, t \in [0,T]$
	\begin{align*}
			D^{H_k}_{r}\Big(\int_{0}^{t} f_j(\boldsymbol{Y}(s)) \mathrm{d}s\Big)
			&=\int_{0}^{t} D^{H_k}_{r}(f_j(\boldsymbol{Y}(s))) \mathrm{d}s\cr 
			&=\int_{0}^{t} D^{H_k}_{r}(Y_1(s))\frac{\partial f_j}{\partial y_1}(\boldsymbol{Y}(s))\mathrm{d}s
			+\int_{0}^{t} D^{H_k}_{r}(Y_2(s))\frac{\partial f_j}{\partial y_2}(\boldsymbol{Y}(s))\mathrm{d}s,
		\end{align*}
		and
		\begin{align*}
			D^{H_k}_{r}\Big(\int_{0}^{t} \sigma_{j1} \mathrm{d}^{\diamond} B^{H_1}_1(s)
			+\int_{0}^{t} \sigma_{j2} \mathrm{d}^{\diamond} B^{H_2}_2(s)\Big)
			&=\left\{\begin{matrix} 
				H_1(2H_1-1)\int_{0}^{t}\sigma_{j1}\vert r-s\vert^{2H_1-2}\mathrm{d}s,&k=1, \\  
				H_2(2H_2-1)\int_{0}^{t}\sigma_{j2}\vert r-s\vert^{2H_2-2}\mathrm{d}s,&k= 2, 
			\end{matrix}\right. \cr
	\end{align*}
	hold for $j,k=1,2$.
\end{lemma}

\section{Derivatives for the explicit expression of Malliavin derivatives}\label{sec-app-2}
The subsequent step involves deriving the explicit expression for the Malliavin derivatives $D^{H_k}_{t}(Y_j(t))$, where $j,k=1,2$. For a fixed $r$, based on Eq.~(\ref{eq-app-13}), and by applying Lemma \ref{malian-to-fwis} together with the initial condition $D^{H_k}_{r}(\boldsymbol{Y}(0))=\boldsymbol{0}$ for $k=1,2$, we consider two distinct cases\\
\textbf{Case 1 $k=1$}
	\begin{align}\label{Malz-1}
		D^{H_1}_{r}(Y_j(t))&= D^{H_1}_{r}\Big(\int_{0}^{t} f_j(\boldsymbol{Y}(s)) \mathrm{d}s\Big)
		+D^{H_1}_{r}\Big(\int_{0}^{t} \sigma_{j1} \mathrm{d}^{\diamond} B^{H_1}_1(s)
		+\int_{0}^{t} \sigma_{j2} \mathrm{d}^{\diamond} B^{H_2}_2(s)\Big)\cr
		&=\int_{0}^{t}D^{H_1}_{r}(Y_1(s))\frac{\partial f_j}{\partial y_1}(\boldsymbol{Y}(s))\mathrm{d}s
		+\int_{0}^{t}D^{H_1}_{r}(Y_2(s))\frac{\partial f_j}{\partial y_2}(\boldsymbol{Y}(s))\mathrm{d}s\cr
		&\quad+H_1(2H_1-1)\int_{0}^{t}\sigma_{j1}\vert r-s\vert^{2H_1-2}\mathrm{d}s,
	\end{align}
	\textbf{Case 2 $k=2$}
	\begin{align}\label{Malz-2}
		D^{H_2}_{r}(Y_j(t))&= D^{H_2}_{r}\Big(\int_{0}^{t} f_j(\boldsymbol{Y}(s)) \mathrm{d}s\Big)
		+D^{H_2}_{r}\Big(\int_{0}^{t} \sigma_{j1} \mathrm{d}^{\diamond} B^{H_1}_1(s)
		+\int_{0}^{t} \sigma_{j2} \mathrm{d}^{\diamond} B^{H_2}_2(s)\Big)\cr
		&=\int_{0}^{t}D^{H_2}_{r}(Y_1(s))\frac{\partial f_j}{\partial y_1}(\boldsymbol{Y}(s))\mathrm{d}s
		+\int_{0}^{t}D^{H_2}_{r}(Y_2(s))\frac{\partial f_j}{\partial y_2}(\boldsymbol{Y}(s))\mathrm{d}s\cr
		&\quad+H_2(2H_2-1)\int_{0}^{t}\sigma_{j2}\vert r-s\vert^{2H_2-2}\mathrm{d}s.
\end{align}
Eqs.~(\ref{Malz-1}) and (\ref{Malz-2}) can be written in matrix form as follows
\begin{align}\label{Malz-11}
		\mathrm{d}D_{r}^{H_1} (\boldsymbol{Y}(t))=\nabla\boldsymbol{\hat f}(\boldsymbol{Y}(t))D_{r}^{H_1} (\boldsymbol{Y}(t))\mathrm{d}t+
		H_1(2H_1-1)(\sigma_{11},\sigma_{21})^{T}\vert t-r \vert ^{2H_1-2}\mathrm{d}t,
	\end{align}
	and
	\begin{align}\label{Malz-21}
		\mathrm{d}D_{r}^{H_2} (\boldsymbol{Y}(t))=\nabla\boldsymbol{\hat f}(\boldsymbol{Y}(t))D_{r}^{H_2} (\boldsymbol{Y}(t))\mathrm{d}t+
		H_2(2H_2-1)(\sigma_{12},\sigma_{22})^{T}\vert t-r \vert ^{2H_2-2}\mathrm{d}t,
\end{align}
where the Jacobian matrix $\nabla\boldsymbol{\hat f}$ is
\begin{align*}
	\nabla\boldsymbol{\hat f}=\begin{pmatrix}
		\dfrac{\partial f_1}{\partial y_1} & \dfrac{\partial f_1}{\partial y_2} \\
		\dfrac{\partial f_2}{\partial y_1} & \dfrac{\partial f_2}{\partial y_2}
	\end{pmatrix}.
\end{align*}

Because $\nabla\boldsymbol{\hat f}$ is a continuous matrix-valued function of time $t$, it can be shown that the aforementioned Eqs.~(\ref{Malz-11}) and (\ref{Malz-21}) can be interpreted as systems of time-varying linear ordinary differential equations, whose solutions can be expressed in the following form
\begin{align}\label{Malz-13}
	D_{r}^{H_1} (\boldsymbol{Y}(t))=H_1(2H_1-1)\int_0^t\Psi(t;s)
	(\sigma_{11},\sigma_{21})^{T}\vert r-s \vert ^{2H_1-2}\mathrm{d}s,
\end{align}
and
\begin{align}\label{Malz-23}
	D_{r}^{H_2} (\boldsymbol{Y}(t))=H_2(2H_2-1)\int_0^t\Psi(t;s)
	(\sigma_{12},\sigma_{22})^{T}\vert r-s \vert ^{2H_2-2}\mathrm{d}s,
\end{align}
where the state-transition matrix $\Psi(t;s)$ satisfies
\begin{align*}
	\frac{\mathrm{d}\Psi(t;s)}{\mathrm{d}t}&=\nabla\boldsymbol{\hat f}(\boldsymbol{Y}(t))\Psi(t;s),\quad\Psi(s,s)=I.
\end{align*}

The Malliavin derivatives $ D^{H_k}_{t}Y_j(t), j,k = 1, 2 $, are obtained by substituting $ t $ for $ r $ in the expressions provided by Eqs.~(\ref{Malz-13}) and (\ref{Malz-23}). \qed 

\section{Derivation for the memFPK equation (Theorem \ref{them-1})}\label{sec-app-3}
Applying the fractional It\^o's formula (see, Lemma \ref{mul-itoformula}) to the It\^o-type SDE \eqref{eq-app-13} for an arbitrary function $F(\boldsymbol{Y}(t))$, one has
\begin{align}\label{itoful-fwis}
	\mathrm{d}F(\boldsymbol{Y}(t))
	=&\sum_{i=1}^{2}f_i(\boldsymbol{Y}(t))\frac{\partial F}{\partial y_i}(\boldsymbol{Y}(t))\mathrm{d}t
	+\sum_{i=1}^{2}\sum_{k=1}^{2}\sigma_{ik}\frac{\partial F}{\partial y_i}(\boldsymbol{Y}(t))\mathrm{d}^{\diamond}B_k^{H_k}(t)\cr
	&+\sum_{i=1}^{2}\sum_{j=1}^{2}\sum_{k=1}^{2}\sigma_{ik}D^{H_k}_{t}(Y_j(t))\frac{\partial^2 F}{\partial y_i\partial y_j}(\boldsymbol{Y}(t))\mathrm{d}t.
\end{align}

As can be seen from Eqs.~(\ref{Malz-13}) and (\ref{Malz-23}), the explicit expression of the Malliavin derivatives is calculated in the \ref{sec-app-2}. Then, taking the expectation of Eq.~\eqref{itoful-fwis} and noting that
\begin{align*}
	\mathbb{E}\Big[\sum_{i=1}^{2}\sum_{k=1}^{2}\sigma_{ik}\frac{\partial F}{\partial y_i}(\boldsymbol{Y}(t))\mathrm{d}^{\diamond}B_k^{H_k}(t)\Big]=0,
\end{align*}
by the zero-mean property of the FWIS integral \cite{ducan2000stochastic}, we obtain
\begin{align}\label{exp-itoful-fwis}
	\frac{\mathrm{d}}{\mathrm{d}t}\mathbb{E}[F(\boldsymbol{Y}(t))]=&
	\mathbb{E}\Big[\sum_{i=1}^{2}f_i(\boldsymbol{Y}(t))\frac{\partial F}{\partial y_i}(\boldsymbol{Y}(t))\Big]
	+\mathbb{E}\Big[\sum_{i=1}^{2}\sum_{j=1}^{2}\sum_{k=1}^{2}\sigma_{ik}D^{H_k}_{t}(Y_j(t))\frac{\partial^2 F}{\partial y_i\partial y_j}(\boldsymbol{Y}(t))\Big].
\end{align}

According to the definition of expectation, Eq.~(\ref{exp-itoful-fwis}) is equivalent to
\begin{align}\label{exp-itoful-fwis-1}
	\frac{\mathrm{d} }{\mathrm{d} t}\int_{\mathbb{R}^2}F(\boldsymbol{y})p(\boldsymbol{y},t)\mathrm{d}\boldsymbol{y}=&\int_{\mathbb{R}^2}\Big\{\sum_{i=1}^{2}f_i(\boldsymbol{y})\frac{\partial F}{\partial y_i}(\boldsymbol{y})\cr
	&+\sum_{i=1}^{2}\sum_{j=1}^{2}\sum_{k=1}^{2}\sigma_{ik}\mathbb{E}[D^{H_k}_{t}(Y_j(t))\vert\boldsymbol{Y}(t)=\boldsymbol{y}]\frac{\partial^2 F}{\partial y_i\partial y_j}(\boldsymbol{y})\Big\}p(\boldsymbol{y},t)\mathrm{d}\boldsymbol{y}.
\end{align}

For the sake of simplicity, combining Eqs.~(\ref{Malz-13}) and (\ref{Malz-23}), denoting 
\begin{align}\label{exp-itoful-fwis-4}
	B_{jk}(\boldsymbol{y},t)&=
	\mathbb{E}[D^{H_k}_{t}(Y_j(t))\vert\boldsymbol{Y}(t)=\boldsymbol{y}]\cr
	&=H_k(2H_k-1)\int_0^t\mathbb{E}[\sigma_{1k}\Psi_{j1}(t;s)+\sigma_{2k}\Psi_{j2}(t;s)\vert \boldsymbol{Y}(t)=\boldsymbol{y}]\vert t-s \vert ^{2H_k-2} \mathrm{d}s.
\end{align}
Applying the integration by parts of the right-hand side of Eq.~(\ref{exp-itoful-fwis-1}) and assuming that as $  \vert y_i \vert \to \infty,i=1,2 $, $ p(\boldsymbol{y},t) $ and its derivatives up to the second order with respect to $ \boldsymbol{Y} $ all vanish, Eq.~(\ref{exp-itoful-fwis-1}) can be rewritten as
\begin{align}\label{exp-itoful-fwis-2}	
	\int_{\mathbb{R}^2}F(\boldsymbol{y})\frac{\partial }{\partial t}p(\boldsymbol{y},t)\mathrm{d}\boldsymbol{y}
	=&
	-\int_{\mathbb{R}^2}F(\boldsymbol{y}) \sum_{i=1}^{2}\frac{\partial }{\partial y_i}\big\{f_i(\boldsymbol{y})p(\boldsymbol{y},t)\big\}\mathrm{d}\boldsymbol{y}\cr
	&+\int_{\mathbb{R}^2}F(\boldsymbol{y})\sum_{i=1}^{2}\sum_{j=1}^{2}\frac{\partial^2 }{\partial y_i\partial y_j} \big\{\sum_{k=1}^{2}\sigma_{ik}B_{jk}(\boldsymbol{y},t)p(\boldsymbol{y},t)\big\}\mathrm{d}\boldsymbol{y}.
\end{align}	

Using the fact that $ F(\boldsymbol{Y}(t)) $ is arbitrary, Eq.~(\ref{exp-itoful-fwis-2}) leads to the following form
\begin{align}\label{exp-itoful-fwis-3}	
	\frac{\partial }{\partial t}p(\boldsymbol{y},t)
	=
	-\sum_{i=1}^{2}\frac{\partial }{\partial y_i}\big\{f_i(\boldsymbol{y})p(\boldsymbol{y},t)\big\}
	+\sum_{i=1}^{2}\sum_{j=1}^{2}\frac{\partial^2 }{\partial y_i\partial y_j} \big\{\sum_{k=1}^{2}\sigma_{ik}B_{jk}(\boldsymbol{y},t)p(\boldsymbol{y},t)\big\}.
\end{align}	

Thus, Eq.~(\ref{exp-itoful-fwis-3}) together with Eq.~(\ref{exp-itoful-fwis-4}) lead to the proof of Theorem \ref{them-1}.

\end{sloppypar}

\end{document}